\documentclass{amsart}

\usepackage{amsthm}
\usepackage{amsfonts}
\usepackage{amsmath}
\usepackage{amssymb}
\usepackage{mathrsfs}
\usepackage{fullpage}
\usepackage{stmaryrd}
\usepackage{hyperref}
\usepackage{verbatim}

\theoremstyle{plain}

\theoremstyle{definition}

\theoremstyle{remark}

\newcommand{\Begin}[2]{\begin{#1}\label{#2}}

\newcommand{\bPi}{\mathbf{\Pi}}
\newcommand{\bSigma}{\mathbf{\Sigma}}
\newcommand{\bDelta}{\mathbf{\Delta}}

\newcommand{\bbP}{\mathbb{P}}
\newcommand{\bbQ}{\mathbb{Q}}
\newcommand{\bbR}{\mathbb{R}}
\newcommand{\bbS}{\mathbb{S}}
\newcommand{\bbU}{\mathbb{U}}
\newcommand{\bbO}{\mathbb{O}}
\newcommand{\bbX}{\mathbb{X}}

\newcommand{\CM}{\mathcal{M}}

\newcommand{\forces}{\Vdash}
\newcommand{\analytic}{{\bSigma_1^1}}
\newcommand{\lanalytic}{{\Sigma_1^1}}

\newcommand{\borel}{{\bDelta_1^1}}

\newcommand{\cantorspace}{{{}^\omega 2}}
\newcommand{\bairespace}{{{}^\omega\omega}}
\newcommand{\finBinarySequence}{{{}^{<\omega}2}}

\newcommand{\reals}{\bbR}
\newcommand{\Los}{\L{}o\'s}
\newcommand{\OD}{\mathrm{OD}}
\newcommand{\HOD}{\mathrm{HOD}}

\newcommand{\degrees}{{\mathcal{D}}}

\begin{document}

\title{Ordinal Definability and Combinatorics of Equivalence Relations}

\author{William Chan}
\address{Department of Mathematics, University of North Texas, Denton, TX 76203}
\email{William.Chan@unt.edu}

\begin{abstract}
Assume $\mathsf{ZF + AD^+ + V = L(\mathscr{P}(\reals))}$. Let $E$ be a $\analytic$ equivalence relation coded in $\HOD$. $E$ has an ordinal definable equivalence class without any ordinal definable elements if and only if $\HOD \models E$ is unpinned. 

$\mathsf{ZF + AD^+ + V = L(\mathscr{P}(\reals))}$ proves $E$-class section uniformization when $E$ is a $\analytic$ equivalence relation on $\reals$ which is pinned in every transitive model of $\mathsf{ZFC}$ containing the real which codes $E$: Suppose $R$ is a relation on $\reals$ such that each section $R_x = \{y : (x,y) \in R\}$ is an $E$-class, then there is a function $f : \reals \rightarrow \reals$ such that for all $x \in \reals$, $R(x,f(x))$. 

$\mathsf{ZF + AD}$ proves that $\reals \times \kappa$ is J\'onsson whenever $\kappa$ is an ordinal: For every function $f : [\reals \times \kappa]^{<\omega}_= \rightarrow \reals \times \kappa$, there is an $A \subseteq \reals \times \kappa$ with $A$ in bijection with $\reals \times \kappa$ and $f[[A]^{<\omega}_=] \neq \reals \times \kappa$. 
\end{abstract}

\maketitle\let\thefootnote\relax\footnote{November 12, 2017. The author was supported by NSF grant DMS-1703708.}


\section{Introduction}\label{Introduction}

The questions of concern here are problems of independent interests that appeared during the study of the J\'onsson property for nonwellorderable sets under the axiom of determinacy.

Let $N \in \omega \cup \{\omega\}$ and $X$ be some set. Define $[X]^N_= = \{x \in {}^NX : (\forall i,j < N)(i \neq j \Rightarrow x(i) \neq x(j))\}$ and $[X]^{<\omega}_= = \bigcup_{n \in \omega} [X]^n_=$. Let $\approx$ denote the relation of being in bijection. Define $\mathscr{P}^N(X) = \{Y \subseteq X : Y \approx N\}$ and $\mathscr{P}^{<\omega}(X) = \bigcup_{n \in \omega} \mathscr{P}^n(X)$. 

An $N$-Jonsson function for $X$ is a function $f : [X]^N_= \rightarrow X$ so that for all $Y \subseteq X$ with $Y \approx X$, $f[[Y]^N_=] = X$. A function $f : [X]^{<\omega}_= \rightarrow X$ is a J\'onsson function if and only if for all $Y \subseteq X$ with $Y \approx X$, $f[[Y]^{<\omega}_=] = X$. A set $X$ has the J\'onsson property if and only if there are no J\'onsson functions for $X$. 

The classical study of the J\'onsson property involved wellordered sets. For wellordered sets $X$, J\'onsson functions for $X$ are formulated using $\mathscr{P}^N(X)$ rather than $[X]^N_=$. Under $\mathsf{AC}$, the following results are known: \cite{On-a-Problem-of-B-Jonsson} showed that every infinite set has an $\omega$-J\'onsson function. The existence of such a function is also where Kunen's proof of the Kunen's inconsistency uses $\mathsf{AC}$. The existence of a cardinal with the J\'onsson property implies $0^\sharp$ exists. Results of Erd\H{o}s and Hajnal (see \cite{Some-Weak-Versions-of-Large-Cardinal-Axioms} and \cite{On-a-Problem-of-B-Jonsson}) imply that under $\mathsf{CH}$, $2^{\aleph_0}$ is not J\'onsson. Hence $\reals$ is not J\'onsson under $\mathsf{CH}$. On the other hand, real valued measurable cardinals are J\'onsson (see \cite{Some-Weak-Versions-of-Large-Cardinal-Axioms} Corollary 11.1). Solovay showed it is consistent relative to a measurable cardinal that $2^{\aleph_0}$ is real valued measurable. Hence it is consistent relative to a measurable cardinal that $\reals$ is J\'onsson. 

Using the axiom of determinacy $\mathsf{AD}$, \cite{Infinitary-Combinatorics-and-the-Axiom-of-Determinateness} showed that $\aleph_n$ is J\'onsson for each $n \in \omega$. \cite{Determinacy-and-Jonsson-Cardinals-in-LR} showed that every cardinal $\kappa < \Theta$ is J\'onsson under $\mathsf{ZF + AD + V = L(\reals)}$. In fact, Woodin showed that $\mathsf{ZF + AD^+}$ can prove every cardinal $\kappa < \Theta$ is J\'onsson.

Under $\mathsf{AD}$, there are sets which cannot be wellordered. Some important examples are quotients of $\borel$ equivalence relations such as $=$, $E_0$, $E_1$, $E_2$, and $E_3$ (see Definition \ref{some examples of equivalence relations}). Holshouser and Jackson (see \cite{Partition-Properties-for-Hyperfinite-Quotients} and \cite{Partition-Properties-for-Non-Ordinal-Sets}) showed that $\reals$ has the J\'onsson property and there are no $2$-J\'onsson functions for $\reals \slash E_0$ under $\mathsf{AD}$. \cite{Definable-Combinatorics-Some-Borel-Equivalence-Relations} showed that under $\mathsf{AD}$, there is a $3$-J\'onsson function for $\reals \slash E_0$. Results from \cite{Definable-Combinatorics-Some-Borel-Equivalence-Relations} seem to suggest that $\reals \slash E_1$, $\reals \slash E_2$, and $\reals \slash E_3$ do not have that J\'onsson property, but no J\'onsson functions for these quotients have yet to be constructed.

For the $\borel$ equivalence relations mentioned above, various dichotomy theorems assert the significance of these equivalence relations in the degree structure of $\borel$ equivalence relations under $\borel$ reducibility. The proofs of these dichotomy results give specific combinatorial structures to sets $A$ such that $E \leq_\borel E \upharpoonright A$, when $E$ is one of the $\borel$ equivalence relations above. For example, if $A \subseteq \reals$ is $\analytic$ and $E_0 \leq_\borel E_0 \upharpoonright A$, then $A$ contains an $E_0$-tree (a perfect tree with very specific symmetry conditions; see \cite{Definable-Combinatorics-Some-Borel-Equivalence-Relations} Definition 5.2). Similarly, if $A \subseteq \reals$ is $\analytic$ and $E_2 \leq_\borel E_2 \upharpoonright A$, then $A$ contains an $E_2$-tree (a perfect tree with certain summability conditions; see \cite{Definable-Combinatorics-Some-Borel-Equivalence-Relations} Fact 14.14).

The following describes the techniques from \cite{Definable-Combinatorics-Some-Borel-Equivalence-Relations} for investigating the J\'onsson property for $\reals \slash E_0$: To study functions $f : [\reals \slash E_0]^2_= \rightarrow \reals \slash E_0$, one would like to lift $f$ to a function $F : \reals^2 \rightarrow \reals$ with the property that for all $(x_1,x_2) \in \reals^2$, $[F(x_1, x_2)]_{E_0} = f([x_1]_{E_0}, [x_2]_{E_0})$. Such a function $F$ is called a lift of $f$. Then one tries to produce an $E_0$-tree on which the collapse of $F$ misses elements of $\reals \slash E_0$. On the other hand, using the specific combinatorial structure of $E_0$-trees, one can define a map $F : \reals^3 \rightarrow \reals$ which is $E_0$-invariant and given any real $x$, there is a triple $(x_1,x_2,x_3)$ of $E_0$-unrelated reals so that $F(x_1,x_2,x_3) \ E_0 \ x$. The collapse of $F$ would then be a $3$-J\'onsson map.

As described in the above example, the existence of lifts of functions from $\reals \slash E \rightarrow \reals \slash F$, where $E$ and $F$ are equivalence relations on $\reals$, seems to be useful in the study of functions on quotients. The existence of a lift is an immediate consequence of uniformization. $\mathsf{AD}_\reals$ has full uniformization. Moreover, a lift of a function $f : \reals \slash E \rightarrow \reals \slash F$ requires only uniformization for relations whose sections are $F$-classes. Woodin showed that countable section uniformization holds in $\mathsf{AD}^+$. Thus lifts exist for functions into $\reals \slash E_0$ under $\mathsf{AD}^+$. Moreover for showing that there are no $2$-J\'onsson functions for $\reals \slash E_0$, it suffices to apply comeager uniformization (which holds in just $\mathsf{AD}$) to find a function $F : C \rightarrow \reals$, where $C \subseteq \reals^2$ is comeager, which lifts $f$ on $C$. Such a lift is adequate since the $2$-Mycielski property for $E_0$ shows that there is a set $A$ such that $E_0 \leq_\borel E_0 \upharpoonright A$ and $\{(x_1,x_2) \subseteq A^2 : \neg(x_1 \ E_ 0 \ x_2)\} \subseteq C$. This roughly implies that $F$ lifts $f$ on a set whose quotient by $E_0$ has cardinality $\reals \slash E_0$. However, \cite{Definable-Combinatorics-Some-Borel-Equivalence-Relations} showed that except for $=$ which has the full Mycielski property, a very limited amount of the Mycielski property holds for the other equivalence relations of interest.

Motivated by this question of $E$-class section uniformization, Zapletal asked a related question: Does every ordinal definable $E_2$ equivalence class contain an ordinal definable real, under $\mathsf{ZF + AD + V = L(\reals)}$? He informed the author that the equivalence relation $=^+$, defined on ${}^\omega\reals$ as equality of range, has ordinal definable classes with no ordinal definable elements and that this phenomenon can be viewed as a consequence of the unpinnedness of $=^+$. He asked then whether pinnedness can be used to characterize those $\borel$ equivalence relations with ordinal definable equivalence classes without any ordinal definable elements. 

Under $\mathsf{AD}$, every ordinal definable countable set of reals contains only ordinal definable elements. The proof of this can be found within the proof of Woodin's countable section enumeration under $\mathsf{AD}^+$, which states that for every relation $R$ with countable sections there is a function that takes $x$ to a wellordering of the section $R_x$. The main idea is to consider the canonical wellordering of $R_x$ in $\HOD_S^{L[S,x,z]}$ as $z$ ranges over a Turing cone of reals and $S$ is some set of ordinals from an $\infty$-Borel code for $R$. (See \cite{Ramsey-Ultrafiler-and-Countable-to-One-Uniformation} for the proof.) This implies that under $\mathsf{AD^+}$, every ordinal definable $E$ class contains only ordinal definable elements if $E$ is an equivalence relation with all countable classes defined using only ordinal parameters. $\mathsf{AD}$ is important for these questions since \cite{A-Definable-E0-Class-Containing-No-Definable-Elements} showed that in a forcing extension of the constructible universe $L$, there is an ordinal definable $E_0$ equivalence class with no ordinal definable elements. 

Section \ref{Ordinal Definable Equivalence Classes} will show roughly that in $L(\reals) \models \mathsf{AD}$, if a $\lanalytic$ equivalence relation $E$ has an $\mathrm{OD}$ equivalence class without any $\OD$ elements, then $\HOD$ must think that $E$ is unpinned:
\\*
\\*
\noindent\textbf{Theorem \ref{OD class no OD element implies unpinned}}
\textit{Assume $\mathsf{ZF + AD^+ + V = L(\mathscr{P}(\reals))}$. Let $T$ be a set of ordinals. Let $E$ be an equivalence relation which is $\lanalytic(s)$ for some $s \in \HOD_T$ and let $A$ be an $\OD_T$ $E$-class. If $A$ has no $\OD_T$ elements, then $\HOD_T \models E$ is unpinned.}
\\*
\\*\indent Models of $\mathsf{ZF + AD^+ + V = L(\mathscr{P}(\reals))}$ are considered natural models of $\mathsf{AD^+}$. If $L(\reals) \models \mathsf{AD}$, then $L(\reals)$ satisfies this theory. Woodin, \cite{A-Trichotomy-Theorem-in-Natural} Corollary 3.2, has shown that if $\mathsf{ZF + AD^+ + V = L(\mathscr{P}(\reals))}$ holds, then either there is a set of ordinals $J$ so that $V = L(J,\reals)$ or else $V \models \mathsf{AD_\reals}$.

The proof of this theorem uses the idea of taking ultraproducts of $\HOD_{S}^{L[S,z]}$ (where the Turing degree of $z$ serves as the index and $S$ is a set of ordinals) using Martin's Turing cone measure. This technique appears in Woodin's proof that sets of reals have $\infty$-Borel codes in $L(\reals)$ when $L(\reals) \models \mathsf{AD}$ as exposited in \cite{Proper-Forcing-and-Absoluteness-Comm} Claim 1.6.
\\*
\\*\noindent\textbf{Theorem \ref{unpinned in HOD implies OD class no OD elements}} 
\textit{$(\mathsf{ZF + AD^+})$ Let $E$ be a $\analytic$ equivalence relation defined in $\HOD_R$, where $R$ is some set. Suppose $\HOD_R \models E$ is unpinned. Then there is an $\OD_R$ $E$-class with no $\OD_R$ elements.}
\\*
\\*\indent These two results together give a very succient answer to Zapletal's question in natural models of $\mathsf{AD^+}$: 
\\*
\\*\noindent\textbf{Corollary \ref{od class unpinnedness equiv}}
\textit{Assume $\mathsf{ZF + AD^+ + V = L(\mathscr{P}(\reals))}$. Let $E$ be a $\analytic$ equivalence relation coded in $\HOD$. $E$ has an $\OD$ $E$-class with no $\OD$ elements if and only if $\HOD \models E$ is unpinned.}
\\*
\\*\indent Many important examples of pinned $\borel$ equivalence relations include $=$, $E_0$, $E_1$, $E_2$, smooth, hyperfinite, and hypersmooth equivalence relations. 

Using the previous theorem, one obtains $E$-class section uniformization for equivalence relations satisfying some definable pinnedness condition. This is particular useful when the equivalence relations are provably pinned:
\\*
\\*\noindent\textbf{Theorem \ref{pinned equiv class section uniformization}}
\textit{Assume $\mathsf{ZF + AD^+ + V = L(\mathscr{P}(\reals))}$. If $E$ is a $\analytic$ equivalence relation which is pinned in every transitive model of $\mathsf{ZFC}$ containing the real that codes $E$, then every relation $R$ whose sections are all $E$-classes can be uniformized.}
\\*
\\*\indent As a consequence, every function $f : \reals \slash E \rightarrow \reals \slash F$ has a lift under $\mathsf{AD^+ + V = L(\mathscr{P}(\reals))}$ when $F$ is $=$, $E_0$, $E_1$, $E_2$, smooth, hyperfinite, essentially countable, or hypersmooth.

Section \ref{jonsson property} will study the J\'onsson property of some nonwellorderable sets. Holshouser and Jackson have shown that $\reals \times \kappa$ for any $\kappa < \Theta$ has the J\'onsson property. They use that $\reals$ and all ordinals $\kappa < \Theta$ have the J\'onsson property. A natural question would be whether $\reals \times \kappa$ is J\'onsson for all ordinals $\kappa$. The proof that $\reals$ is J\'onsson has a clear flavor of classical descriptive set theory since it uses comeagerness, continuity, the Mycielski property, and fusions of perfect trees. The proof that ordinals $\kappa < \Theta$ is J\'onsson have a somewhat different flavor. A related question would be whether the J\'onsson property for $\kappa$ is relevant to showing $\reals \times \kappa$ is J\'onsson. Does there exists a more classical proof that $\reals \times \kappa$ is J\'onsson? It will be shown that:
\\*
\\*\textbf{Theorem \ref{R product WO is jonsson}}
\textit{$(\mathsf{ZF + AD})$ For any ordinal $\kappa$, $\reals \times \kappa$ has the J\'onsson property.}
\\*
\\*\indent Whether or not $\kappa$ is J\'onsson does not appear in the proof of the above theorem. This result is proved while investigating the J\'onsson property for wellordered disjoint unions $\bigsqcup_{\alpha < \kappa} \reals \slash E_\alpha$ where each $E_\alpha$ is an equivalence relation with all classes countable and $\reals \slash E_\alpha \approx \reals$. The techniques have a very classical flavor using results about lengths of wellordered sequences of reals, additivity of the meager ideal, comeager uniformization, and fusions of perfect trees. There are also some discussions about the cardinality of $\bigsqcup_{\alpha < \kappa} \reals \slash E_\alpha$. However, it remains open whether $\bigsqcup_{\alpha < \kappa} \reals \slash E_\alpha$ has the J\'onsson property. 

This section concludes by producing a $6$-J\'onsson function for $(\reals \slash E_0) \times \kappa$ for any $\kappa < \Theta$ under $\mathsf{AD}$. This shows that $(\reals \slash E_0) \times \kappa$ for $\kappa < \Theta$ is not J\'onsson under $\mathsf{AD}$.

The author would like to thank Jared Holshouser, Stephen Jackson, Alexander Kechris, Connor Meehan, and Itay Neeman for comments and discussions about the material in this paper. In particular, the author would like to thank Jind\v{r}ich Zapletal for informing the author of the main question.

\section{Ordinal Definable Equivalence Classes}\label{Ordinal Definable Equivalence Classes}

$V$ will denote the universe of set theory in consideration. If $M$ is a model of set theory and $A$ is some concept given by some formula, then $A^M$ will denote the relativization of that formula inside $M$. If a concept $A$ is unrelativized, then it is assumed to mean $A^V$, although it may be written $A^V$ for emphasis. $\bbR$ will denote $\bairespace$, the Baire space, consisting of functions from $\omega$ to $\omega$ with its usual metric. (Although it may sometimes denote $\cantorspace$, the Cantor space.) The elements of $\bbR$ will be called reals.

If $X$ is a set, then $\OD_X$ denotes the class of sets which are ordinal definable using $X$ as a parameter. $\HOD_X$ is the collection of sets which are hereditarily ordinal definable from $X$. $\HOD_X \models \mathsf{ZFC}$ and has a canonical global wellordering definable using $X$.

\Begin{fact}{vopenka result}
(Vop\v{e}nka) Suppose $S$ is a set of ordinals. Let $x \in \reals$. 

In $L[S,x]$, let $\bbP$ denote the forcing of $\OD_S$ subsets of $\reals$ ordered by $\subseteq$. Using the canonical $S$-definable bijection of $\OD_S$ subsets onto $\mathrm{ON}$, let $\bbO_S \in \HOD_S$ be the forcing that results by transferring $\bbP$ onto $\mathrm{ON}$ using this map.

Then there is a $G \in L[S,x]$, which is $\bbO_S$-generic over $\HOD_S$, so that $L[S,x] = \HOD_S[G] = \HOD_S[x]$. 
\end{fact}

\begin{proof}
See \cite{Set-Theory} Theorem 15.46.
\end{proof}

\Begin{definition}{infinity Borel code}
Let $X \subseteq \reals$, $S$ be a set of ordinals, and $\varphi$ be a formula in the language of set theory. $(S,\varphi)$ is an $\infty$-Borel code for $X$ if and only if for all $x \in \reals$, $x \in X \Leftrightarrow L[S,x] \models \varphi(S,x)$. 
\end{definition}

\Begin{definition}{AD+}
(\cite{Axiom-of-Determinacy-Forcing-Axioms} Section 9.1) $\mathsf{AD}^+$ consists of the following:

\noindent (1) $\mathrm{DC}_\bbR$.

\noindent (2) Every $A \subseteq \bbR$ has an $\infty$-Borel code.

\noindent (3) For all $\lambda < \Theta$, $A \subseteq \reals$, and continuous function $\pi : {}^\omega \lambda \rightarrow \reals$,  $\pi^{-1}[A]$ is determined. 

($\lambda$ is given the discrete topology. $\Theta$ is the supremum of the ordinals which are surjective images of $\reals$. Games with moves from $\lambda$ are defined the same way as the more familiar games on $\omega$.)
\end{definition}

\Begin{definition}{pinned concepts}
(\cite{Reducibility-Invariants-in-Higher-Set-Theory}) Let $E$ be an equivalence relation on $\reals$. Let $\bbP$ be a forcing. Let $\tau$ be a $\bbP$-name. 

Let $\tau_\mathrm{left}, \tau_\mathrm{right}$ be the canonical $\bbP \times \bbP$-names with the property that $\tau_\mathrm{left}$ and $\tau_\mathrm{right}$ are evaluated according to $\tau$ using the left and right $\bbP$-generic filters, respectively, coming from a $\bbP\times\bbP$-generic filter.

$\tau$ is an $E$-pinned name if and only if $1_{\bbP \times \bbP} \forces_{\bbP \times \bbP} \tau_\mathrm{left} \ E \ \tau_\mathrm{right}$.

An $E$-pinned name $\tau$ is an $E$-trivial name if and only if there is some $x \in \bbR$ so that $1_\bbP \forces_{\bbP} \tau \ E \ \check x$. 

$E$ is a pinned equivalence relation if and only if all $E$-pinned names are $E$-trivial.
\end{definition}

Pinnedness is more accurately a property of a fixed definition for the equivalence relation $E$ (which is to be used to interpret $E$ in generic extensions). This paper is concern only with $\analytic$ equivalence relations and such equivalence relations are always defined as the projection of certains trees on $\omega \times \omega \times \omega$.

\Begin{definition}{martin measure}
Let $\leq_T$ denote the Turing reducibility relation on $\bairespace$. For $x,y \in \bairespace$, let $x \equiv_T y$ if and only if $x \leq_T y$ and $y \leq_T x$. A Turing degree is a $\equiv_T$ equivalence class. If $x,y \in \bairespace$, then define $[x]_{\equiv_T} \leq_T [y]_{\equiv_T}$ if and only if $x \leq_T y$. 

Let $\degrees$ denote the set of Turing degrees. A Turing cone with base $C \in \degrees$ is the set $\{D \in \degrees : C \leq_T D\}$. Define Martin's measure $\mathcal{U}$ by: for $A \in \mathscr{P}(\degrees)$, $A \in \mathcal{U}$ if and only if $A$ contains a Turing cone.

Under $\mathsf{AD}$, the Martin's measure is a countably complete ultrafilter on $\degrees$.
\end{definition}

\Begin{definition}{od ultrapower}
$(\mathsf{ZF + AD})$ Let $T$ be some set. Let $\mathcal{H}$ be a (usually proper class) function on $\degrees$ which is definable using only $T$ and ordinals as parameters and takes each $X$ to some transitive class. Assume that there is some (usually proper class) function $\mathcal{R}$ definable using only $T$ and ordinals as parameters so that for each $X \in \degrees$, $\mathcal{R}(X)$ is a wellordering of $\mathcal{H}(X)$. 

Let $M^T_{\mathcal{H}, \mathcal{R}}$ denote the collection of $\OD_T$ functions on $\degrees$ taking each $X \in \degrees$ to an element in $\mathcal{H}(X)$. For $F,G \in M^T_{\mathcal{H},\mathcal{R}}$, let $F \sim G$ if and only if $\{X \in \degrees : F(X) = G(X)\} \in \mathcal{U}$. 

Let $\CM^T_{\mathcal{H},\mathcal{R}}$ denote the collection of equivalence classes of $M^T_{\mathcal{H},\mathcal{R}}$ under $\sim$. Define $[F]_\sim \in  [G]_\sim$ if and only if $\{X \in \degrees : F(X) \in G(X)\} \in \mathcal{U}$. 
\end{definition}

\Begin{fact}{properties of od ultraproduct}
$(\mathsf{ZF +AD})$ $\CM_{\mathcal{H},\mathcal{R}}^T$ is a $T$-definable class consisting of $\OD_T$ elements. Using the $T$-definable bijection of $\OD_T$ and $\mathrm{ON}$, $\CM^T_{\mathcal{H},\mathcal{R}}$ is isomorphic to a class inside $\HOD_T$. $\CM_{\mathcal{H},\mathcal{R}}^T$ is well-founded; hence, it can be considered as a transitive structure inside $\HOD_T$. 

The \Los's theorem holds for $\CM_{\mathcal{H},\mathcal{R}}^T$: Suppose $F_0, ..., F_{k - 1} \in M^T_{\mathcal{H},\mathcal{R}}$ and $\varphi$ is a formula of $\{\dot \in\}$, then $\CM^T_{\mathcal{H},\mathcal{R}} \models \varphi([F_0]_\sim, ..., [F_k]_{\sim})$ if and only $\{X \in \degrees : \mathcal{H}(X) \models \varphi(F_0(X), ..., F_{k - 1}(X))\} \in \mathcal{U}$. 

For each $\alpha < \omega_1$, let $c_\alpha : \mathcal{D} \rightarrow \{\alpha\}$ be the constant function taking value $\alpha$. The class $[c_\alpha]_\sim$ represents the ordinal $\alpha$ in $\CM_{\mathcal{H},\mathcal{R}}^T$. 

For each $r \in \reals$ which is $\OD_T$ and belongs to $\mathcal{H}(X)$ for a cone of $X \in \degrees$, define the function $c_r : \degrees \rightarrow \{\emptyset, r\}$ by $c_r(X) = r$ if $r \in \mathcal{H}(X)$ and $c_r(X) = \emptyset$ if otherwise. Then $[c_r]_\sim$ represents $r$ in $\CM_{\mathcal{H},\mathcal{R}}^T$. 
\end{fact}

\begin{proof}
$\CM^T_{\mathcal{H},\mathcal{R}}$ is a structure in $\OD_T$ since $M^T_{\mathcal{H},\mathcal{R}} \subseteq \mathrm{OD}_T$. Note the $\in$ relation of $\CM^T_{\mathcal{H},\mathcal{R}}$ is definable from $T$. Using the the canonical bijection of $\OD_T$ and $\mathrm{ON}$, one can transfer $\CM_{\mathcal{H},\mathcal{R}}^T$ and its $\in$-relation onto $\mathrm{ON}$. This new isomorphic structure consists entirely of ordinals and hence elements of $\HOD_T$. 

Let $F \in M^T_{\mathcal{H},\mathcal{R}}$. Suppose $[F]_\sim$ is not wellfounded. There is some set $X \subseteq \{[G]_\sim : [G]_\sim \in [F]_\sim\}$ without an $\in^{\CM_{\mathcal{H},\mathcal{R}}^{T}}$-minimal element. Let $L(0)$ be the $\OD_T$-least function $G$ so that $[G]_\sim \in X$. Suppose $L(n)$ has been defined. Let $L(n + 1)$ be the $\OD_T$-least function $G$ so that $[G]_\sim \in X$ and $[G]_\sim \in [L(n)]_\sim$. Let $A_n = \{x \in \degrees : L(n + 1)(x) \in L(n)(x)\}$. Each $A_n \in \mathcal{U}$. Since $\mathcal{U}$ is countably complete, $\bigcap_{n \in \omega} A_n \neq \emptyset$. Let $x \in \bigcap_{n \in \omega} A_n$. Then $\langle L(n)(x) : n \in \omega\rangle$ is an $\in$-decreasing sequence in $V$. Contradiction. $\CM^T_{\mathcal{H},\mathcal{R}}$ is well-founded. Using the Mostowski collapse, one may consider $\CM^T_{\mathcal{H},\mathcal{R}}$ as a transitive structure inside of $\HOD_T$. 

The proof of \Los's theorem is by induction on formula complexity: The result holds for the atomic formulas by definition. Assume the result holds for $\varphi$ and $\psi$, then the result holds for $\neg\varphi$ and $\varphi \wedge \psi$ by the usual arguments. (Note the case involving $\neg$ requires that $\mathcal{U}$ is an ultrafilter.) Suppose the result has been shown for $\varphi$. If $\CM_{\mathcal{H},\mathcal{R}}^T \models (\exists x)\varphi(x, [F_0]_\sim, ..., [F_{k - 1}]_\sim)$, then there exists some $G \in M_{\mathcal{H},\mathcal{R}}^T$ so that $\CM_{\mathcal{H},\mathcal{R}}^T \models \varphi([G]_\sim, [F_0]_\sim, ..., [F_{k - 1}]_\sim)$. Using the induction hypothesis, $\{X \in \degrees : \mathcal{H}(X) \models (\exists x)\varphi(x, F_0(X), ..., F_{k - 1}(X))\} \in \mathcal{U}$. Suppose $\{X \in \degrees : (\exists x)\varphi(x, F_0(X), ..., F_{k - 1}(X))\} \in \mathcal{U}$. Define $G$ on $\degrees$ by letting $G(X)$ be the $\mathcal{R}(X)$-least element $z$ of $\mathcal{H}(X)$ such that $\mathcal{H}(X) \models \varphi(z, F_0(X), ..., F_{k - 1}(X))$ if such an element exists and $\emptyset$ otherwise. $G$ is $\OD_T$ and so belongs to $M^T_{\mathcal{H},\mathcal{R}}$. By the induction hypothesis, $\CM^T_{\mathcal{H},\mathcal{R}} \models \varphi([G]_\sim, [F_0]_\sim, ..., [F_{k - 1}]_\sim)$. Therefore, $\CM^T_{\mathcal{H},\mathcal{R}} \models (\exists x)\varphi(x, [F_0]_\sim, ..., [F_{k - 1}]_\sim)$. This completes the sketch of \Los's theorem.

Suppose $[F]_\sim \in [c_\alpha]_\sim$. Let $A = \{X \in \degrees : F(X) \in \alpha\}$. $A \in \mathcal{U}$. Let $A_\beta= \{X \in \degrees : F(X) = \beta\}$. $\mathcal{A} = \bigcup_{\beta < \alpha} A_\beta$. Since $\mathcal{U}$ is countably complete and $\alpha$ is countable, there is some $\beta < \alpha$ so that $A_\beta \in \mathcal{U}$. Then $c_\beta \sim F$. This shows that $[c_\alpha]_\sim$ represents $\alpha$ in $\CM^T_{\mathcal{H},\mathcal{R}}$ when $\alpha < \omega_1$. 
\end{proof}

\Begin{fact}{location infinity borel code}
(Woodin, \cite{A-Trichotomy-Theorem-in-Natural} Theorem 3.4) Assume $\mathsf{ZF + AD^+ + V = L(\mathscr{P}(\reals))}$. Let $T$ be a set of ordinals. A set $X \subseteq \reals$ which is $\OD_T$ has an $\infty$-Borel code $(S,\varphi)$ which is $\OD_T$.
\end{fact}

\Begin{fact}{HOD is constructible from ordinal}
(Woodin, \cite{A-Trichotomy-Theorem-in-Natural} Theorem 2.18) Assume $\mathsf{ZF + AD^+ + V = L(\mathscr{P}(\reals))}$. Let $T$ be a set of ordinals. There is some set of ordinals $\bbX$ so that $\HOD_T = L[\bbX]$. (Note that $\bbX$ is $\OD_T$.)
\end{fact}

In the case of $L(\reals)$ and $T = \emptyset$, the set $\bbX$ can be taken to be $\bbP^\omega$ which is the direct limit indexed by $n \in \omega$ of Vop\v{e}nka forcing on $\reals^n$. This follows from Woodins result that $L(\reals)$ is a symmetric collapse extension of its $\HOD$. One can find an exposition of this result in \cite{Proper-Forcing-and-Absoluteness-Comm}.

\Begin{fact}{wellfoundedness ordinal ultraproduct}
(Woodin, \cite{A-Trichotomy-Theorem-in-Natural} Section 2.2) Assuming $\mathsf{ZF + AD^+}$, $\prod_{X \in \degrees}\mathrm{ON} \slash \mathcal{U}$ is wellfounded.
\end{fact}

Assume $\mathsf{AD^+}$, the wellfoundedness of $\CM^T_{\mathcal{H},\mathcal{R}}$ can also be proved from Fact \ref{wellfoundedness ordinal ultraproduct}. For the question of Zapletal, one will need to form an ultraproduct of the form $\CM^T_{\mathcal{H},\mathcal{R}}$ so that all the reals of $\HOD$ belong to this ultraproduct.

\Begin{fact}{hod reals in specific ultraproduct}
Assume $\mathsf{ZF + AD^+ + V = L(\mathscr{P}(\reals))}$. Let $T$ be a set of ordinals. Let $\bbX$ be a set of ordinals as given by Fact \ref{HOD is constructible from ordinal}, so that $\HOD_T = L[\bbX]$. For each $X \in \degrees$, let $\mathcal{H}(X) = \HOD_\bbX^{L[\bbX,X]}$ and $\mathcal{R}(X)$ be the canonical wellordering of $\HOD_\bbX^{L[\bbX,X]}$. Then $\CM^{\bbX}_{\mathcal{H},\mathcal{R}}$ is wellfounded, $\CM^\bbX_{\mathcal{H},\mathcal{R}} \subseteq \HOD_T$, and $\reals^{\HOD_T} \subseteq \mathcal{M}^T_{\mathcal{R},\mathcal{H}}$. 
\end{fact}

\begin{proof}
Note that $\bbX$ is $\OD_T$. Observe that for all $X \in \degrees$, $\HOD_T = L[\bbX] \subseteq \HOD_\bbX^{L[\bbX,X]}$. So if $r \in \HOD_T$, then $r \in \HOD_{\bbX}^{L[\bbX,X]}$. The function $c_r$ is $\OD_\bbX$ and belongs to $M^\bbX_{\mathcal{H},\mathcal{R}}$. This result now follows from Fact \ref{properties of od ultraproduct}.
\end{proof}

\Begin{theorem}{OD class no OD element implies unpinned}
Assume $\mathsf{ZF + AD^+ + V = L(\mathscr{P}(\reals))}$. Let $T$ be a set of ordinals. Let $E$ be an equivalence relation which is $\lanalytic(s)$ for some $s \in \HOD_T$ and let $A$ be an $\OD_T$ $E$-class. If $A$ has no $\OD_T$ elements, then $\HOD_T \models E$ is unpinned.
\end{theorem}

\begin{proof}
For simplicity, let $T = \emptyset$. By Fact \ref{HOD is constructible from ordinal}, let $\bbX$ be a set of ordinals so that $\HOD = L[\bbX]$. By Fact \ref{location infinity borel code}, $A$ has an $\infty$-Borel code in $\HOD = L[\bbX]$. Modifying $\bbX$ by including an ordinal if necessary, one may as well assume that there is some formula $\varphi$ so that $(\bbX,\varphi)$ forms an $\infty$-Borel code for $A$. 

Recall that $E$ is $\lanalytic(s)$ means there is some $s$-recursive tree $T$ on $\omega \times \omega \times \omega$ so that $x \ E \ y$ if and only if $L[s,x,y] \models T^{x,y}$ is illfounded, where $T^{x,y} = \{u : (x \upharpoonright |u|, y \upharpoonright |u|, u) \in T\}$. In this way, $E$ is $\infty$-Borel with a code that is a subset of $\omega$. 

Suppose $y \geq_T x$ for some $x \in A$. By Fact \ref{vopenka result}, there is some $\bbO_{\bbX}^{L[\bbX,y]}$-name $\tau \in \HOD^{L[\bbX,y]}_{\bbX}$ and some $\bbO_{\bbX}^{L[\bbX,y]}$-generic over $\HOD_{\bbX}^{L[\bbX,y]}$ filter $G \in L[\bbX,y]$ so that $\tau[G] = x$ and $L[\bbX,y] = \HOD_\bbX^{L[\bbX,y]}[G]$. Since $V \models L[\bbX,x] \models \varphi(\bbX,x)$, $L[\bbX,y] \models L[\bbX,x] \models \varphi(\bbX,x)$. Since $L[\bbX,y] = \HOD_{\bbX}^{L[\bbX,y]}[G]$, one has $\HOD_{\bbX}^{L[\bbX,y]}[G] \models L[\bbX,x] \models \varphi(\bbX,x)$. There is some $q \in \bbO^{L[\bbX,y]}_{\bbX}$ so that $\HOD_{\bbX}^{L[\bbX,y]} \models q \forces_{\bbO_{\bbX}} L[\check \bbX,\tau] \models \varphi(\check \bbX,\tau)$. Let $q_y$ and $\tau_y$ be the $\HOD^{L[\bbX,y]}_\bbX$-least such $q$ and $\tau$ with the above properties. In order to satisfy the technical requirement of using the largest condition of the forcing in the definition of pinnedness, let $\bbU_y = \{p \in \bbO^{L[\bbX,y]}_\bbX : p \leq_{\bbO^{L[\bbX,y]}_\bbX} q_y\}$, $\leq_{\bbU_y} = \leq_{\bbO^{L[\bbX,y]}_\bbX} \upharpoonright \bbU_y$, and $1_{\bbU_y} = q_y$. If $y$ does not Turing compute any element of $A$, then one can just let $\bbU_y$ and $\tau_y$ be $\emptyset$.

If $x \equiv_T y$, $\HOD^{L[\bbX,x]}_\bbX = \HOD^{L[\bbX,y]}_\bbX$ and their canonical global wellorderings are the same. This shows that $\bbU_x = \bbU_y$ and $\tau_x = \tau_y$. If $X \in \degrees$ and $x \in X$, let $\HOD^{L[\bbX,X]}_\bbX = \HOD^{L[\bbX,x]}_\bbX$, $\bbU_X = \bbU_x$, and $\tau_X = \tau_x$. For $X \in \degrees$, let $\mathcal{H}(X) = \HOD^{L[\bbX,X]}_\bbX$ and $R(X)$ be the canonical global wellordering of $\HOD_\bbX^{L[\bbX,X]}$. For $X \in \degrees$, define $\Phi_\bbU(X) = \bbU_X$ and $\Phi_\tau(X) = \tau_X$. Let $\CM = \CM^\bbX_{\mathcal{H},\mathcal{R}}$. Note that $\Phi_\bbU,\Phi_\tau \in M_{\mathcal{H},\mathcal{R}}^\bbX$. Let $\bbU = [\Phi_\bbU]_\sim$ and $\tau = [\Phi_\tau]_\sim$. Let $c_\bbX$ be the constant function taking value $\bbX$. Note that $c_\bbX \in M^\bbX_{\mathcal{H},\mathcal{R}}$. Let $\bbX^\infty = [c_\bbX]_\sim$. As in Fact \ref{properties of od ultraproduct}, $\CM$ will be identified as a transitive class in $\HOD^V$. Thus $\bbU$, $\tau$, and $\bbX^\infty$ belong to $\HOD^V$.

By \Los's theorem, $\CM$ is a model of $\mathsf{ZFC}$, $\bbU$ is some forcing, $\tau$ is some $\bbU$-name adding a real, $\bbX^\infty$ is a set of ordinals, and $\CM \models 1_{\bbU} \forces_\bbU L[\bbX^\infty,\tau] \models \varphi(\bbX^\infty,\tau)$. 

Claim 1: 
$$\CM \models 1_{\bbU\times\bbU} \forces_{\bbU \times \bbU} (\forall x)(\forall y)((L[\bbX^\infty,x] \models \varphi(\bbX^\infty,x) \wedge L[\bbX^\infty,y]\models \varphi(\bbX^\infty,y)) \Rightarrow x \ E \ y)$$

(Note that the ultraproduct moves $\bbX$ to $\bbX^\infty$. However, $E$ as a $\lanalytic(s)$ equivalence relation has the real $s$ as its $\infty$-Borel code. The constant function $c_s$ taking value $s$ belongs to $M^\bbX_{\mathcal{H},\mathcal{R}}$. In $\CM$, $[c_s]_\sim$ represents $s$. That is, $s$ is not moved by the ultraproduct. Hence it is appropriate to continue to denote $E$ by $E$ in $\CM$ as it is still the same $\analytic$ equivalence relation.)

To see the claim: Fix some $z \in A$. By \Los's theorem, it suffices to prove that for all $r \geq_T z$:
$$\HOD_\bbX^{L[\bbX,r]} \models 1_{\bbU_r \times \bbU_r} \forces_{\bbU_r \times \bbU_r} (\forall x)(\forall y)((L[\bbX,x] \models \varphi(\bbX,x) \wedge L[\bbX,y] \models \varphi(\bbX,y)) \Rightarrow x \ E \ y)$$
Fix some $(p,q) \in \bbU_r \times \bbU_r$. Since $L[\bbX,r] \models \mathsf{AC}$ and $V \models \mathsf{AD}$, $\omega_1^V$ is inaccessible in $\HOD^{L[\bbX,r]}_\bbX$. Hence $\bbU_r \times \bbU_r$ and its power set in $\HOD^{L[\bbX,r]}_\bbX$ are countable in $V$. There exists $G \times H \in V$ which is $\bbU_r\times\bbU_r$-generic over $\HOD_\bbX^{L[\bbX,r]}$. Since $G \times H \in V$, all sets of $\HOD^{L[\bbX,r]}_\bbX[G\times H]$ belong to $V$. Let $x$ and $y$ be reals of $\HOD_\bbX^{L[\bbX,r]}[G\times H]$ so that $\HOD_\bbX^{L[\bbX,r]}[G \times H] \models L[\bbX,x] \models \varphi(\bbX,x) \wedge L[\bbX,y] \models \varphi(\bbX,y)$. Then $V \models L[\bbX,x] \models \varphi(\bbX,x) \wedge L[\bbX,y] \models \varphi(\bbX,y)$. Since $(\bbX,\varphi)$ is an $\infty$-Borel code for $A$ in $V$, $x \in A$ and $y \in A$. Since $A$ is an $E$-class, $x \ E \ y$. By Mostowski absolutenss, $\HOD_\bbX^{L[\bbX,r]}[G\times H] \models x \ E \ y$. This shows that $\HOD^{L[\bbX,r]}_\bbX[G\times H]$ satisfies the formula behind the above forcing relation. Since $G \times H$ is generic, there is some $(p',q') \leq_{\bbU_r \times \bbU_r} (p,q)$ so that in $\HOD^{L[\bbX,r]}_\bbX$, $(p',q')$ forces that formula. Since $(p,q)$ was arbitrary, this establishes the claim.

Claim 2: 
$$\CM \models 1_{\bbU \times \bbU} \forces_{\bbU \times \bbU} (\forall x)(\forall y)((L[\bbX^\infty,x] \models \varphi(\bbX^\infty,x) \wedge x \ E \ y) \Rightarrow L[\bbX^\infty,y] \models \varphi(\bbX^\infty,y))$$
The proof essentially uses the same idea as Claim 1.

Now to show that $\bbU$ and $\tau$ witness that $E$ is unpinned in $\HOD^V$: 

First to show that $\tau$ is an $E$-pinned name in $\HOD^V$: Let $G \times H$ be any $\bbU \times \bbU$-generic filter over $\HOD^V$. Since $\CM \subseteq \HOD^V$, if $G$ and $H$ are generic over $\HOD^V$, then $G$ and $H$ are generic over $\CM$. By the forcing theorem, $\CM[G] \models L[\bbX^\infty,\tau[G]] \models \varphi(\bbX^\infty,\tau[G])$ and $\CM[G] \models L[\bbX^\infty,\tau[H]] \models \varphi(\bbX^\infty,\tau[H])$. By Claim 1, $\CM[G \times H] \models \tau[G] \ E \ \tau[H]$. By Mostowski absoluteness, $\HOD^V[G \times H] \models \tau[G] \ E \ \tau[H]$. Since $G \times H$ was arbitrary, $\HOD_\bbX^V \models 1_{\bbU \times \bbU} \forces_{\bbU\times\bbU} \tau_\mathrm{left} \ E \ \tau_\mathrm{right}$. This shows that $\tau$ is an $E$-pinned $\bbU$-name in $\HOD^V$. 

Finally, to show that $\tau$ is not $E$-trivial: Suppose there is some $x \in \HOD^V$ so that $\HOD^V \models 1_\bbU \forces_\bbU \tau \ E \ \check x$. Let $G \subseteq \bbU$ be a $\bbU$-generic over $\HOD^V$ filter. Then $\HOD^V[G] \models \tau[G] \ E \ x$. By Mostowski absoluteness, $\CM[G] \models \tau[G] \ E \ x$. $G$ is also generic over $\CM$. By the forcing theorem, $\CM[G] \models L[\bbX^\infty,\tau[G]] \models \varphi(\bbX^\infty,\tau[G])$. Since $x \in \HOD^V$, Fact \ref{hod reals in specific ultraproduct} and Fact \ref{properties of od ultraproduct} imply that $[c_x]_\sim$ represents $x$ in $\CM$. By Claim 2 applied in $\CM[G \times H]$ where $H$ is any $\bbU$-generic filter over $\CM[G]$, $\CM[G] \models L[\bbX^\infty,x] \models \varphi(\bbX^\infty,x)$. Thus $\CM \models L[[c_\bbX]_\sim, [c_x]_\sim] \models \varphi([c_\bbX]_\sim,[c_x]_\sim)$. By \Los's theorem, for a Turing cone of $X$'s (such that $x \in \HOD^{L[\bbX,X]}_\bbX$), $\HOD^{L[\bbX,X]}_\bbX \models L[\bbX,x] \models \varphi(\bbX,x)$. This implies $V \models L[\bbX,x] \models \varphi(\bbX,x)$. $V \models x \in A$ since $(\bbX,\varphi)$ is the $\infty$-Borel code for $A$ in $V$. This contradicts the assumption that $A$ has no $\OD$ elements.

This completes the proof.
\end{proof}

\Begin{theorem}{unpinned in HOD implies OD class no OD elements}
$(\mathsf{ZF + AD})$ Let $E$ be a $\analytic$ equivalence relation defined in $\HOD_R$, where $R$ is some set. Suppose $\HOD_R \models E$ is unpinned. Then there is an $\OD_R$ $E$-class with no $\OD_R$ elements.
\end{theorem}

\begin{proof}
Since $\HOD_R \models E$ is unpinned, there exists some forcing $\bbP \in \HOD_R$ and $\bbP$-name $\sigma \in \HOD_R$ so that within $\HOD_R$, $\bbP$ and $\sigma$ witness that $E$ is not pinned. 

Inside $\HOD_R$ (which models $\mathsf{AC}$), let $N$ be an elementary substructure of some large enough rank initial segment of $\HOD_R$ with the property that (1) $N$ contains the code for $E$, (2) $\bbR \subseteq N$, (3) $\bbP,\sigma \in N$, and (4) $N$ has cardinality $|\bbR|$. Let $M$ be the Mostowski collapse of $N$. Let $\bbQ$ and $\tau$ be the image of $\bbP$ and $\sigma$ under the Mostowski collapse map. As $E$ is $\analytic$, the code for $E$ is a tree on $\omega \times \omega \times \omega$ whose projection is $E$. So a code for $E$ is merely a subset of $\omega$. Hence the Mostowski collapse map does not move the code for $E$. Note that $|M|^V = |\bbR^{\HOD_R^{V}}|^V = \aleph_0$ since $\mathsf{AD}$ holds. Hence there are generics for $\bbQ$ over $M$ that lie in $V$. 

Suppose $G$ and $H$ are two generic filters for $\bbQ$ over $M$ which belong to $V$. Since $M[G]$ and $M[H]$ are countable in $V$, one can construct a generic filter $J \in V$ so that $G \times J$ and $H \times J$ are generic filters for $\bbQ \times \bbQ$. By elementarity, $M \models \tau$ is $E$-pinned. Thus $M[G \times J] \models \tau[G] \ E \ \tau[J]$ and $M[H \times J] \models \tau[H] \ E \ \tau[J]$. By Mostowski absoluteness, $\tau[G] \ E \ \tau[J]$ and $\tau[H] \ E \ \tau[J]$ holds in $V$. Since $E$ is an equivalence relation, $\tau[G] \ E \ \tau[H]$. This shows that whenever $G$ and $H$ are $\bbQ$-generic filters over $M$ that belong to $V$ (but may not be mutually generic), $\tau[G] \ E \ \tau[H]$. 

$M \models \tau$ is not $E$-trivial by elementarity. Since $\bbR^{\HOD_R} \subseteq M$, for any $G \subseteq \bbQ$ which is $\bbQ$-generic over $M$ and any $x \in \bbR^{\HOD_R}$, $M[G] \models \neg(\tau[G] \ E \ x)$. By absoluteness, if $G \in V$, then $\neg(\tau[G] \ E \ x)$. 

In $V$, let $A$ be the set of $x \in \bbR$ so that there exists some $G \subseteq \bbQ$ which is $\bbQ$-generic over $M$ and $x \ E \ \tau[G]$. Since $\bbQ,\tau \in M$ and $M \in \HOD_R$, $A$ is $\OD_R$. By the discussion of the above two paragraphs, $A$ is a single $E$-class and has no elements of $\OD_R$. 

Note that the only consequence of $\mathsf{AD}$ that is used is that there is no uncountable wellordered set of reals.
\end{proof}

The following answers the question of Zapletal. 

\Begin{corollary}{od class unpinnedness equiv}
Assume $\mathsf{ZF + AD^+ + V = L(\mathscr{P}(\reals))}$. Let $E$ be a $\analytic$ equivalence relation coded in $\HOD$. $E$ has an $\OD$ $E$-class with no $\OD$ elements if and only if $\HOD \models E$ is unpinned.
\end{corollary}

The rest of this section will give some examples. 

\Begin{definition}{some examples of equivalence relations}
The following are some important $\borel$ equivalence relations.

Let $=$ denote the identity equivalence relation on $\reals$. 

Let $=^+$ denote the Friedman-Stanley jump of $=$ which is defined on ${}^\omega\reals$ by $x =^+ y$ if and only if $\{x(n) : n \in \omega\} = \{y(n) : n \in \omega\}$. ($=^+$ is equality of range.)

Let $E_0$ be the equivalence relation on $\bbR$ (or $\cantorspace$) defined by $x \ E_0 \ y$ if and only if $(\exists k)(\forall n \geq k)(x(n) = y(n))$.

Let $E_1$ be the equivalence relation on ${}^\omega\bbR$ defined by $x \ E_1\ y$ if and only $(\exists k)(\forall n \geq k)(x(n) = y(n))$. 

Let $E_2$ be the equivalence relation on $\cantorspace$ defined by $x \ E_2 \ y$ if and only if $\sum\{\frac{1}{n} : n \in x \triangle y\} < \infty$, where $\triangle$ denotes the symmetric difference operation.
\end{definition}

\Begin{fact}{example pinned and unpinned}
The equivalence relations $=$, $E_0$, $E_1$, and $E_2$ are pinned $\borel$ equivalence relations. Every $\borel$ equivalence relation with countable classes is pinned. Every smooth, hyperfinite, essentially countable, or hypersmooth equivalence relation is pinned.

The equivalence relation $=^+$ is unpinned.
\end{fact}

\begin{proof}
See \cite{Borel-Equivalence-Relations} Chapter 11.

The Solovay product lemma states: Let $\bbP$ and $\bbQ$ be two forcings. Suppose $G \times H$ is $\bbP \times \bbQ$-generic over $V$. Then $V[G] \cap V[H] = V$.

From the Solovay product lemma, it follows that $=$, $E_0$, and $E_1$ are pinned equivalence relations.

If $E \leq_\borel F$ and $F$ is pinned, then $E$ is also pinned. This implies that smooth, hyperfinite, and hypersmooth equivalence relations are pinned.

\cite{Borel-Equivalence-Relations} Theorem 17.1.3 (iii) states that $\borel$ equivalence relations with all classes $\bSigma_3^0$ are pinned. This implies that $E_2$ and every $\borel$ equivalence relation with countable classes are pinned. Therefore, essentially countable equivalence relations are pinned.

Let $\bbQ = \mathrm{Coll}(\omega, \bbR)$. Let $\tau$ be the name for the generic surjection of $\omega$ onto $\bbR$. $\bbQ$ and $\tau$ witness that $=^+$ is unpinned since if $\tau$ was forced to be $=^+$ related to a ground model element, then $\bbR$ would be countable in the ground model.
\end{proof}

\Begin{example}{F2 OD class no OD element}
The proof above that $=^+$ is unpinned can be used to produce an $\OD$ $=^+$-class with no $\OD$ elements assuming $(\mathscr{P}(\bbR))^{\HOD}$ is countable.

Let $\bbQ = \mathrm{Coll}(\omega,\bbR)$ and $\tau$ be the generic surjection of $\omega$ onto $\bbR$ as defined inside of $\HOD$. (Note that $\tau$ is an $=^+$-pinned name.) By the assumption, there exists $\bbQ$-generics over $\HOD$ in $V$. Let $A$ be the collection of $x \in {}^\omega\bbR$ such that there exist some $G \subseteq \bbQ$ which is $\bbQ$-generic over $\HOD$ and $x =^+ \tau[G]$. $A$ is an $\OD$ $=^+$ equivalence class. $A$ cannot contain any $\OD$ elements for otherwise $\HOD$ would think $\bbR^\HOD$ is countable.
\end{example}

\section{Equivalence Class Section Uniformization and Lifting}\label{equiv class section uniform and lifting}

\Begin{theorem}{pinned equiv class section uniformization}
Assume $\mathsf{ZF + AD^+ + V = L(\mathscr{P}(\reals))}$. Let $T$ be a set of ordinals. Let $E$ be a $\analytic$ equivalence relation coded in $\HOD_T$. Suppose $E$ is pinned in $\HOD_{T,x}$ for all $x \in \reals$. Let $R \subseteq \reals \times \reals$ be $\OD_T$ and have the property that for all $x \in \reals$, $R_x = \{y : R(x,y)\}$ is an $E$-class. Then there is a function $F : \reals \rightarrow \reals$ which is $\OD_T$ and uniformizes $R$: that is, for all $x \in \reals$, $R(x,F(x))$. 

If $E$ is a $\analytic$ equivalence relation which is pinned in every transitive model of $\mathsf{ZFC}$ containing the real that codes $E$, then every relation $R$ whose sections are $E$-classes can be uniformized. (For example, $E$ could be any of the pinned equivalence relations from Fact \ref{example pinned and unpinned}.)
\end{theorem}

\begin{proof}
Under these assumptions , for each $x \in \reals$, $R_x$ is an $\OD_{T,x}$ $E$-class. Since $\HOD_{T,x} \models E$ is unpinned, Theorem \ref{OD class no OD element implies unpinned} implies that $R_x$ must have an $\OD_{T,x}$ element. For each $x \in \reals$, let $F(x)$ be the least element of $\HOD_{T,x}$ under the canonical global wellordering of $\HOD_{T,x}$ which belongs to $R_x$. $F$ is an $\OD_{T}$ uniformization of $R$. 

For the second statement, under $\mathrm{AD}^+$, any such relation $R$ has an $\infty$-Borel code $(S,\varphi)$. By modifying $S$ if necessary, one may assume that $\HOD_S$ contains a code for $E$ as a $\analytic$ set. By the hypothesis, $E$ is pinned in every $\HOD_{S,x}$, where $x \in \reals$. The second statement follows from the first statement.
\end{proof}

\cite{Reducibility-Invariants-in-Higher-Set-Theory} has shown that if $E$ is a $\borel$ equivalence relation coded in some transitive model $M$ and $N$ is some transitive model with $M \subseteq N$, then $E$ is pinned in $M$ if and only if $E$ is pinned in $N $. Therefore, in the first statement of Theorem \ref{pinned equiv class section uniformization}, it suffices just to have $\HOD_{T} \models E$ is pinned, when $E$ is a $\borel$ equivalence relation. 

However \cite{Reducibility-Invariants-in-Higher-Set-Theory} also shows that, in general, pinnedness for $\analytic$ equivalence relation is not absolute by producing a pinned $\analytic$ equivalence relation in $L$ which is unpinned in a forcing extension of $L$. However, in the present situation, one is concerned with models of the form $\HOD_T^V$ and $\HOD_{T,x}^V$ where $V$ is a model of determinacy. Possible more can be said in such settings. This suggests the following question.

\Begin{question}{pinned in HOD condition replace}
In the first statement of Theorem \ref{pinned equiv class section uniformization}, can the condition that $E$ is pinned in $\HOD_{T,x}$ for all $x \in \reals$ be replace by just $E$ is pinned in $\HOD_{T}$ when $E$ is a $\analytic$ equivalence relation coded in $\HOD_T$?
\end{question}

Regardless, most natural examples are $\borel$. Moreover, for most of the natural examples, pinnedness is provable in $\mathsf{ZFC}$. 

\Begin{definition}{lifting maps on quotients}
Let $E$ be an equivalence relation on some set $X$. Let $F$ be an equivalence relation on some set $Y$. Let $n \in \omega$. Let $f : (X \slash E)^n \rightarrow (Y \slash F)$ be some function. A function $F : X^n \rightarrow Y$ is a lift of $f$ if and only if for all $x_0, ..., x_{n - 1} \in X$, $[F(x_0, ..., x_{n - 1})]_F = f([x_0]_E, ..., [x_{n - 1}]_E)$. 
\end{definition}

\Begin{corollary}{pinned equivalence relation lifting}
Assume $\mathsf{ZF + AD^+ + V = L(\mathscr{P}(\reals))}$. Suppose $E$ is an equivalence relation on $\bbR$. Suppose $F$ is a $\analytic$ equivalence relation on $\bbR$ which is pinned in every transitive models of $\mathsf{ZFC}$ containing the real that codes $F$. For all $n \in \omega$, every function $f : (\reals \slash E)^n \rightarrow (\reals \slash F)$ has a lift.

In particular, this lifting property holds when $F$ is $E_0$, $E_1$, $E_2$, smooth, hyperfinite, essentially countable, or hypersmooth.
\end{corollary}

\begin{proof}
Define the relation $R(x_0, ..., x_{n - 1}, y)$ if and only if $y \in f([x_0]_E, ..., [x_{n - 1}]_E)$. For each $(x_0, ..., x_{n - 1}) \in \reals^n$, $R_{(x_0, ..., x_{n - 1})} = f([x_0]_E, ..., [x_{n - 1}]_E)$, which is an $F$-class. By assumption, $F$ is pinned in every model of $\mathsf{ZFC}$ containing the real that codes $F$. Theorem \ref{pinned equiv class section uniformization} implies that $R$ has a uniformizing function $G$. $G$ is a lift of $f$.
\end{proof}

\Begin{example}{AD+ unpinned no uniformization lifting}
Under $\mathsf{ZF + AD_\bbR}$, every relation can be uniformized. Hence, $E$-class section uniformization and lifting for $E$ holds for every equivalence relation $E$ on $\reals$. However $\mathsf{ZF + AD^+}$ is not able to prove $E$-class section uniformization when $E$ is an unpinned equivalence relation. The following is an example.

Assume $\mathsf{ZF + AD + V = L(\reals)}$. 

Define $R(x,y)$ if and only if $y$ is not $\OD_x$. $R$ has no uniformizing function: Suppose $f : \bbR \rightarrow \bbR$ uniformized $R$. Since $V = L(\reals)$, every set of reals is ordinal definable from some real. Thus $f$ is $\OD_z$ for some $z \in \reals$. Hence $f(z)$ is $\OD_z$. However, $R(z,f(z))$ implies that $f(z)$ is not $\OD_z$. Contradiction.

Define $S(x,y)$ if and only if $\{y_n : n \in \omega\} = \bbR^{\HOD_x}$, where $y_n \in \reals$ denotes the $n^\text{th}$ section of $y$ under some coding of pairs of integers by integers. If $S(x,y)$, then $y \notin \OD_x$ for otherwise $\bbR^{\HOD_x}$ would be countable in $\HOD_x$. Since $S \subseteq R$ and $R$ has no uniformization, $S$ also has no uniformization. 
\end{example}

Every instance of $F$-class section uniformization gives a lift of a function from $f : \bbR \rightarrow (\reals \slash F)$. Therefore, failure of $F$-class section uniformization is a failure of lifting for $F$. However, the more interesting instance of the lifting property involving function of the form $f : (\reals \slash F) \rightarrow (\reals \slash F)$. This suggest the following question which may yield more information on the relationship between lifting and $F$-class section uniformization.

\Begin{question}{lift =+ to =+}
Under $\mathsf{ZF + AD + V = L(\reals)}$, is there some function $f : ({}^\omega\reals \slash =^+) \rightarrow ({}^\omega\reals \slash =^+)$ which does not have a lift?
\end{question}

\section{J\'onsson Property}\label{jonsson property}

\Begin{definition}{jonsson property}
Let $X$ be a set and $n \in \omega$. Let $E$ be an equivalance relation on $X$. Let $[X]^n_E = \{(x_0, ..., x_{n - 1}) \in {}^nX : (\forall i < n)(\forall j < n)(i \neq j \Rightarrow \neg(x_i \ E \ x_j))\}$. Let $[X]^{<\omega}_E = \bigcup_{n \in \omega} [X]^n_E$. 

A set $X$ has the J\'onsson property if and only if for all functions $f : [X]^{<\omega}_= \rightarrow X$, there is some $Y \subseteq X$ with $Y \approx X$ and $f[[Y]^{<\omega}_=] \neq X$. (The symbol $\approx$ is the relation of being in bijection.)

For $n < \omega$, an $n$-J\'onsson function for $X$ is a map $f : [X]^n_= \rightarrow X$ so that for all $Y \subseteq X$ with $Y \approx X$, $f[[X]^n_=] = X$. 
\end{definition}

\Begin{fact}{jonsson properties of some sets}
Under $\mathsf{ZF + AD}$, 

(\cite{Partition-Properties-for-Non-Ordinal-Sets} and \cite{Partition-Properties-for-Hyperfinite-Quotients}) $\reals$ has the J\'onsson property. 

(\cite{Definable-Combinatorics-Some-Borel-Equivalence-Relations}) There is a $3$-J\'onsson function for $\reals \slash E_0$. Hence $\reals \slash E_0$ does not have the J\'onsson property. 
\end{fact}

For the rest of this section, $\bbR$ will refer to $\cantorspace$, the set of infinite binary sequences. 

\Begin{definition}{Sacks forcing}
A nonempty subset $p$ of $\finBinarySequence$ is a tree if and only if for all $s \in p$ and $t \subseteq s$, $t \in p$. A tree $p$ is a perfect tree if and only if for all $s \in p$, there is a $t \supseteq s$ so that $t\hat{\ }0, t\hat{\ }1 \in p$. 

Let $\bbS$ be the set of all perfect trees. Let $\leq_\bbS = \subseteq$. 

Let $p \in \bbS$. A node $s \in p$ is a split node if and only if $s\hat{\ }0, s\hat{\ }1 \in p$. A node $s \in p$ is a split of $p$ if and only if $s \upharpoonright (|s| - 1)$ is a split node of $p$. For $n \in \omega$, $s$ is an $n$-split of $p$ if and only if $s$ is a $\subseteq$-minimal element of $p$ with exactly $n$-many proper initial segments which are split nodes of $p$. 

Let $\mathrm{split}^n(p)$ denote the set of $n$-splits of $p$. Note that $|\mathrm{split}^n(p)| = 2^n$ and $\mathrm{split}^0(p) = \{\emptyset\}$. 

If $p,q \in \bbS$, define $p \leq_\bbS^n q$ if and only if $p \leq_\bbS q$ and $\mathrm{split}^n(p) = \mathrm{split}^n(q)$. 

If $p \in \bbS$ and $s \in p$, then define $p_s = \{t \in p : t \subseteq s \vee s \subseteq t\}$. 

Let $p \in \bbS$. Let $\Lambda$ be defined as follows:

\noindent (i) $\Lambda(p,\emptyset) = \emptyset$. 

\noindent (ii) Suppose $\Lambda(p,s)$ has been defined for all $s \in {}^n2$. Fix an $s \in {}^n2$ and $i \in 2$. Let $t \supseteq \Lambda(p,s)$ be the minimal split node of $p$ extending $\Lambda(p,s)$. Let $\Lambda(p, s\hat{\ }i) = t\hat{\ }i$. 

Let $\Xi(p,s) = p_{\Lambda(p,s)}$. 
\end{definition}

\Begin{fact}{fusion lemma}
A fusion sequence is a sequence $\langle p_n : n \in \omega\rangle$ in $\bbS$ so that for all $n \in \omega$, $p_{n + 1} \leq_\bbS^n p_n$. Let $p_\omega = \bigcap_{n \in \omega} p_n$. Then $p_\omega \in \bbS$ and is called the fusion of the above fusion sequence.
\end{fact}

\Begin{fact}{functions unequal tuples miss something}
Suppose $p \in \bbS$. Let $\langle r_n : n \in \omega\rangle$ be a sequence of positive integers. Let $\langle f_n : n \in \omega\rangle$ be a sequence such that for all $n \in \omega$, $f_n : [[p]]^{r_n}_= \rightarrow {}^\omega \bbR$ is a continuous function. Then there is some $q \leq_\bbS p$ and $z \in \bbR$ so that for all $m,n \in \omega$ and $y \in f_n[[[q]]^{r_n}_=]$, $z \neq y(m)$. 
\end{fact}

\begin{proof}
Let $B : \omega \rightarrow \omega \times \omega$ be a surjection with the property that the inverse image of any $(e,g)$ is infinite.

Objects $\langle z_n : n \in \omega \rangle$ and $\langle q_n : n \in \omega\rangle$ will be built with the following properties.

\noindent (I) For each $n \in \omega$, $z_n \in \finBinarySequence$ and $z_n \subsetneq z_{n + 1}$. For each $n \in \omega$, $q_n \in \bbS$, $q_n \leq_\bbS p$, and $q_{n + 1} \leq_\bbS^n q_n$.

\noindent (II) For each $n \in \omega$, suppose $B(n) = (e,g)$. Then for each sequence $(\sigma_1, ..., \sigma_{r_e})$ of pairwise distinct strings in ${}^n2$, there is some $\tau \in \finBinarySequence$ so that for all $y$ with 
$$y \in f_e[[\Xi(q_{n + 1}, \sigma_1)] \times ... \times [\Xi(q_{n + 1},\sigma_{r_e})]]$$
$y(g) \in N_\tau$ and $z_{n + 1}$ and $\tau$ are incompatible. 

Suppose these objects can be constructed. Then $\langle q_n : n \in \omega\rangle$ forms a fusion sequence. By Fact \ref{fusion lemma}, $q = \bigcap_{n \in \omega} q_n$ is a perfect tree. Let $z = \bigcup_{n \in \omega} z_n$. Let $e,g \in \omega$. Suppose $(x_1, ..., x_{r_e}) \in [[q]]^{r_e}_=$. By the assumption on $B$, there is some $n$ large enough so that $B(n) = (e,g)$ and there are pairwise distinct strings $\sigma_1, ..., \sigma_{r_e} \in {}^n2$ with $\Lambda(q,\sigma_1) \subset x_1$, ..., $\Lambda(q,\sigma_{r_e}) \subset x_{r_e}$. Then by (II), $z_{n + 1}$ is not an initial segment of $y(g)$. Hence $y(g) \neq z$. 

It remains to construct these objects.

Let $z_0 = \emptyset$ and $q_0 = p$. 

Suppose $q_n$ and $z_n$ have been constructed. Suppose that $B(n) = (e,g)$. Enumerate all the $r_e$-tuples of distincts strings in ${}^n2$ as $(\sigma_1^0, ..., \sigma_{r_e}^0)$, ..., $(\sigma_1^M, ..., \sigma_{r_e}^M)$ for some $M \in \omega$. 

Let $s_0 = q_n$. Let $\ell_0 = z_n$. Suppose $s_k$ and $\ell_k$ have been defined for some fixed $k \leq M$. For each $1 \leq i \leq r_e$, let $c_i = \sigma_i^k\hat{\ }\bar{0}$. Let $d_i = \bigcup_{n < \omega} \Lambda(r_k, c_i\upharpoonright n)$. By the continuity of $f_e$ on $[[p]]^{r_e}_=$, there is some $N > n$ so that for all 
$$y \in f_e[[\Xi(s_k, c_1\upharpoonright N)] \times ... \times [\Xi(s_k, c_{r_e}\upharpoonright N)]],$$
$f_e(d_1, ..., d_{r_e})(g) \upharpoonright |\ell_k + 1| \subseteq y(g)$. Define $\ell_{k + 1} = \ell_{k} \hat{\ } (1 - f_e(d_1, ...,d_{r_e})(g)(|\ell_k|))$, that is $\ell_{k + 1}$ extends $\ell_k$ by one using the opposite of the value of $f_e(d_1, ..., d_{r_e})(g)(|\ell_k|)$. Let $s_{k + 1} \leq_{\bbS}^n s_k$ be such that for all $\sigma \in {}^n 2$, if $\sigma = \sigma^k_i$ for some $1 \leq i \leq r_e$, then $\Xi(s_{k + 1},\sigma) = \Xi(s_k, c_i\upharpoonright N)$ and if $\sigma$ is otherwise, then $\Xi(s_{k + 1},\sigma) = \Xi(s_k,\sigma)$. 

Finally, let $q_{n + 1} = s_{M + 1}$ and $z_{n + 1} = \ell_{M + 1}$. This completes the construction.
\end{proof}

\Begin{fact}{union of meager sets}
Let $\delta$ be an ordinal. Let $\langle A_\alpha : \alpha < \delta\rangle$ be a sequence of meager subsets of $\reals$. Define a prewellordering on $\bigcup_{\alpha < \delta}A_\alpha$ by $x \preceq y$ if and only if the least ordinal $\xi$ such that $x \in A_\xi$ is less or equal to the least ordinal $\xi$ such that $y \in A_\xi$. Assume that $\preceq$ as a subset of $\reals \times \reals$ has the Baire property. Then $\bigcup_{\alpha < \delta} A_\alpha$ is meager.

$(\mathsf{ZF + AD})$ Every wellordered union of meager sets is meager.
\end{fact}

\begin{proof}
See \cite{Measure-and-Category-in-Effective}. The second statement follows from the fact that every subset of $\reals \times \reals$ has the Baire property under $\mathsf{AD}$. 
\end{proof}

\Begin{fact}{mycielski theorem}
(Mycielski) Suppose $\langle C_n : n \in \omega\rangle$ is a sequence so that each $C_n$ is a comeager subset of $\reals^n$. Then there is a perfect tree $p$ so that for all $n \in \omega$, $[[p]]_=^n \subseteq C_n$. 
\end{fact}

\Begin{fact}{comeager uniformization}
$(\mathsf{ZF + AD})$ (Comeager uniformization) Let $R \subseteq \reals \times \reals$ be a relation. Then there is a comeager set $C \subseteq \reals$ and a function $f : C \rightarrow \reals$ so that for all $x \in C$, $R(x,f(x))$. 
\end{fact}

\Begin{fact}{countable equiv quotient perfect set}
$(\mathsf{ZF + AD})$ Let $E$ be an equivalence relation on $\reals$ with all classes countable and $\reals \slash E \approx \reals$. Let $p$ be perfect tree. Then $[p] \slash E \approx \reals$.
\end{fact}

\begin{proof}
Note that $[p] \slash E$ injects into $\reals \slash E$ by inclusion. Composing with the bijection then shows that $[p] \slash E$ injects into $\reals$. Let $\Phi : \reals \slash E \rightarrow \reals$ be a bijection. Since $E$ has only countable classes and countable unions of countable sets are countable under $\mathsf{AD}$, $[p] \slash E$ is an uncountable set. Hence $\Phi[[p] \slash E]$ is an uncountable subset of $\reals$. By the perfect set property, there is some perfect tree $q$ so that $[q] \subseteq \Phi[[p] \slash E]$. $\Phi^{-1}$ injects $[q]$ into $[p] \slash E$. Hence $\reals$ injects into $[p] \slash E$. By Cantor-Schr\"oder-Bernstein, $[p] \slash E \approx \reals$. 
\end{proof}

\Begin{fact}{no PWO uncountable set countable PWO class}
$(\mathsf{ZF + AD})$ Let $A \subseteq \reals$. Let $\preceq$ be a prewellordering on $A$. For each $x \in A$, let $[x]_{\preceq} = \{y : x \preceq y \wedge y \preceq x\}$. If for all $x \in A$, $[x]_{\preceq}$ is countable, then $A$ is countable. 
\end{fact}

\begin{proof}
If $A$ is not countable, then by the perfect set property, there is some perfect tree $p$ so that $[p] \subseteq A$. Using the notation from Definition \ref{Sacks forcing}, define $x \sqsubseteq y$ if and only if $\bigcup_{n \in \omega} \Lambda(p,x \upharpoonright n) \preceq \bigcup_{n \in \omega} \Lambda(p,y \upharpoonright n)$. Then $\sqsubseteq$ is a prewellordering on $\reals$ so that for each $x \in \reals$, $[x]_{\sqsubseteq}$ is countable. Let $\beta$ be the length of $\sqsubseteq$. For each $\alpha < \beta$, let $A_\alpha$ be the prewellordering class of $\sqsubseteq$ with rank $\alpha$. $\bigcup_{\alpha < \beta} A_\alpha = \reals$ and each $A_\alpha$ is countable (and hence meager). This is not possible by Fact \ref{union of meager sets}.
\end{proof}

\Begin{question}{question disjoint union smooth jonsson}
(Holshouser-Jackson) $(\mathsf{ZF + AD})$ Let $\kappa$ be an ordinal. Let $\langle E_\alpha : \alpha < \kappa\rangle$ be a sequence of equivalence relations on $\reals$ with all classes countable so that $\reals \slash E_\alpha \approx \reals$. Does the disjoint union $\bigsqcup_{\alpha < \kappa} \reals \slash E_\alpha$ have the J\'onsson property?
\end{question}

Note that one is not given a sequence of bijections $\langle \Phi_\alpha : \alpha < \kappa\rangle$ witnessing $\reals \slash E_\alpha \approx \reals$. With such a sequence of bijections, one can construct a bijection witnessing $\bigsqcup_{\alpha < \kappa} \reals \slash E_\alpha \approx \reals \times \kappa$. In this case, Theorem \ref{R product WO is jonsson} below would imply $\bigsqcup_{\alpha < \kappa} \reals \slash E_\alpha$ has the J\'onsson property. The following is an interesting question.

\Begin{question}{disjoint union smooth quotient isomorphic real product ordinal}
(Holshouser-Jackson) $(\mathsf{ZF + AD})$ Let $\kappa$ be an ordinal. Let $\langle E_\alpha : \alpha < \kappa\rangle$ be a sequence of equivalence relations on $\reals$ with all classes countable so that $\reals \slash E_\alpha \approx \reals$. Is $\bigsqcup_{\alpha < \kappa} \reals \slash E_\alpha \approx \reals \times \kappa$?
\end{question}

The following theorem gives some information concerning the J\'onsson property.

\Begin{theorem}{partial jonsson of dijoint union smooth quotient}
$(\mathsf{ZF + AD})$ Let $\kappa$ be an ordinal. Let $\langle E_\alpha : \alpha < \kappa\rangle$ be a sequence of equivalence relation on $\reals$ with all classes countable. Let $f : [\bigsqcup_{\alpha<\kappa} \reals \slash E_\alpha]^{<\omega}_= \rightarrow \bigsqcup_{\alpha < \kappa} \reals \slash E_\alpha$. Then there is some perfect tree $p$ so that $f[[\bigsqcup_{\alpha < \kappa} [p]\slash E_\alpha]^{<\omega}_=] \neq \bigsqcup_{\alpha < \kappa} \reals \slash E_\alpha$. 
\end{theorem}

\begin{proof}
Let $E$ be the equivalence relation on $\reals \times \kappa$ defined by: $(x,\alpha) \ E \ (y, \beta)$ if and only if $\alpha = \beta$ and $x \ E_\alpha \ y$. Then $\bigsqcup_{\alpha < \kappa} \reals \slash E_\alpha$ is in bijection with the quotient $(\reals \times \kappa) \slash E$. In the following $f$ will be considered as a function taking values in $(\reals \times \kappa) \slash E$. 

Let $X$ be the collection of surjections $\sigma : \{1, ..., n\} \rightarrow \{1, ..., m\}$ where $1 \leq m \leq n$ are integers. For all $\sigma \in X$, let $n(\sigma) = n$ and $m(\sigma) = m$, i.e. $n(\sigma)$ and $m(\sigma)$ indicate the domain and range of $\sigma$, respectively.

For each $\sigma \in X$, define $A^\sigma \subseteq \reals^{m(\sigma)} \times \reals$ by
$$(x_1, ..., x_{m(\sigma)}, y) \in A^\sigma \Leftrightarrow (\exists \alpha_1) ...(\exists \alpha_{n(\sigma)})(\exists \beta)([(y,\beta)]_E = f([(x_{\sigma(1)},\alpha_1)]_E, ..., [(x_{\sigma(n(\sigma))},\alpha_{n(\sigma)})]_E))$$

In the following, fix a wellordering of $n(\sigma)$-tuples of ordinals. For each $(x_1, ..., x_{m(\sigma)})$, the elements of $A^\sigma_{(x_1, ..., x_{m(\sigma)})}$ can be prewellordered as follows: $y_0 \sqsubseteq y_1$ if and only if the least $(\alpha_1, ..., \alpha_{n(\sigma)})$ such that there exists (a unique) $\beta$ with 
$$(y_1,\beta) \in f([(x_{\sigma(1)},\alpha_1)]_E, ..., [(x_{\sigma(n(\sigma))},\alpha_{n(\sigma)})]_E)$$
is less than or equal to the least $(\alpha_1, ..., \alpha_{n(\sigma)})$ such that there exists (a unique) $\beta$ with 
$$(y_2,\beta) \in f([(x_{\sigma(1)},\alpha_1)]_E, ..., [(x_{\sigma(n(\sigma))},\alpha_{n(\sigma)})]_E).$$ 
Let $y \in A^\sigma_{(x_1, ..., x_{m(\sigma)})}$. Let $(\alpha_1, ..., \alpha_{n(\sigma)})$ be the least $n(\sigma)$-tuple of ordinals such that for some (unique) $\beta$, $(y,\beta) \in f([(x_{\sigma(1)},\alpha_1)]_E, ..., [(x_{\sigma(n(\sigma))},\alpha_{n(\sigma)})]_E)$. Then $[y]_{\sqsubseteq} \subseteq \pi_1[f([(x_{\sigma(1)},\alpha_1)]_E, ..., [(x_{\sigma(n(\sigma))},\alpha_{n(\sigma)})]_E)]$, where $\pi_1 : \reals \times \kappa \rightarrow \reals$ is the projection onto the first coordinate. Note $f([(x_{\sigma(1)}, \alpha_1), ..., [(x_{\sigma(n(\sigma))},\alpha_{n(\sigma)})]_E)$ is contained inside of $\reals \times \{\beta\}$ for some $\beta$. Since $E_\beta$ is an equivalence relation with countable classes, $\pi_1[f([(x_{\sigma(1)}, \alpha_1), ..., [(x_{\sigma(n(\sigma))},\alpha_{n(\sigma)})]_E)]$ is countable. It has been shown that each $\sqsubseteq$-prewellordering class is countable. By Fact \ref{no PWO uncountable set countable PWO class}, $A^\sigma_{(x_1, ..., x_{m(\sigma)})}$ is countable for all $(x_1, ..., x_{m(\sigma)})$. 

Fix $\sigma \in X$. Let $R^\sigma \subseteq \reals^{m(\sigma)} \times {}^\omega\reals$ be defined by $(x,h) \in R^\sigma$ if and only if $h : \omega \rightarrow A^\sigma_x$ is a surjection. For all $x \in \reals^{\sigma(m)}$, $R^\sigma_x \neq \emptyset$ since $A^\sigma_x$ is countable. By comeager uniformization (Fact \ref{comeager uniformization}), there is some comeager set $C^\sigma \subseteq \reals^{m(\sigma)}$ and some function $H^\sigma : C^\sigma \rightarrow {}^\omega\reals$ so that $R^\sigma(x,H^\sigma(x))$ for all $x \in C^\sigma$. Using $\mathsf{AD}$, one may assume that $H^\sigma$ is continuous on $C^\sigma$ by choosing a smaller comeager set if necessary.

By the results of Mycielski, Fact \ref{mycielski theorem}, there is some perfect tree $p$ so that $[[p]]^{m(\sigma)}_= \subseteq C^\sigma$ for all $\sigma \in X$. Note that for all $\sigma \in X$, $H^\sigma \upharpoonright [[p]]^{m(\sigma)}_=$ is a continuous function with the property that for all $x \in [[p]]^{m(\sigma)}_=$, $H^\sigma(x) \in {}^\omega \reals$ enumerates $A^\sigma_x$. By Fact \ref{functions unequal tuples miss something}, there is some $q \leq_\bbS p$ and some $z \in \reals$ so that for all $\sigma \in X$, $j \in \omega$, and $(x_1, ..., x_{m(\sigma)}) \in [[q]]^{m(\sigma)}_=$, $H^\sigma(x)(j) \neq z$. 

Now suppose $(r_1,\alpha_1), ..., (r_n,\alpha_n) \in [q] \times \kappa$ are such that $([(r_1,\alpha_1)]_E, ..., [(r_n,\alpha_n)]_E) \in [\bigsqcup_{\alpha < \kappa} [q] \slash E_\alpha]^n_=$. There is some $m \leq n$, $(x_1, ..., x_m) \in [[q]]^m_=$, and surjection $\sigma : \{1, ..., n\} \rightarrow \{1, ..., m\}$ so that $(r_1, ..., r_n) = (x_{\sigma(1)}, ..., x_{\sigma(n)})$. Then $z \notin A^\sigma_{(x_1, ..., x_m)}$ implies that $(z,\beta) \notin f([(r_1,\alpha_1)]_E, ..., [(r_n,\alpha_n)]_E)$ for all $\beta < \kappa$. 

This shows that $f[[\bigsqcup_{\alpha < \kappa} [q] \slash E_\alpha]^{<\omega}_=] \neq \bigsqcup_{\alpha < \kappa} \reals \slash E_\alpha$. 
\end{proof}

Let $p$ be the perfect tree given by Theorem \ref{partial jonsson of dijoint union smooth quotient}. Assume that each $\reals \slash E_\alpha \approx \reals$. By Fact \ref{countable equiv quotient perfect set}, each $[p] \slash E_\alpha \cong \reals$. If $\bigsqcup_{\alpha < \kappa} [p] \slash E_\alpha \approx \bigsqcup_{\alpha < \kappa} \reals \slash E_\alpha$, then Theorem \ref{partial jonsson of dijoint union smooth quotient} would imply $\bigsqcup_{\alpha < \kappa} \reals \slash E_\alpha$ has the J\'onsson property. This suggests the following natural question.

\Begin{question}{quotient by perfect tree isomorphism}
$(\mathsf{ZF + AD})$ Let $\kappa$ be an ordinal. Let $\langle E_\alpha : \alpha < \kappa\rangle$ be a sequence of equivalence relations on $\reals$ with all classes countable and $\reals \slash E_\alpha \approx \reals$ for each $\alpha < \kappa$. Let $p$ be a perfect tree. Is $\bigsqcup_{\alpha < \kappa} \reals \slash E_\alpha \approx \bigsqcup_{\alpha < \kappa} [p] \slash E_\alpha$? 
\end{question}

When all the $E_\alpha$'s are the identity equivalence relation, $=$, then one can exhibit the desired bijection. This gives the following result. 

\Begin{theorem}{R product WO is jonsson}
$(\mathsf{ZF + AD})$ For any ordinal $\kappa$, $\reals \times \kappa$ has the J\'onsson property.
\end{theorem}

\begin{proof}
Let $\langle E_\alpha : \alpha < \omega\rangle$ be a sequence where each $E_\alpha$ is the identity equivalence relation, $=$, on $\reals$. Note that $\bigsqcup_{\alpha < \kappa} \reals \slash E_\alpha \approx \reals \times \kappa$. Apply Theorem \ref{partial jonsson of dijoint union smooth quotient} to this sequence. For any perfect tree $p$, $\bigsqcup_{\alpha < \kappa} [p]\slash E_\alpha \approx \bigsqcup_{\alpha < \kappa} [p] \approx \reals \times \kappa$. 
\end{proof}

Many of the results above are trivial if the sequence $\langle E_\alpha : \alpha < \kappa\rangle$ is accompanied by a sequence $\langle \Phi_\alpha : \alpha < \kappa\rangle$ where each $\Phi_\alpha : \reals \slash E_\alpha \rightarrow \reals$ is a bijection. A natural question would be to construct an example $\langle E_\alpha : \alpha < \kappa\rangle$ such that for each $\alpha < \kappa$, $\reals \slash E_\alpha \approx \reals$ but there does not exists a sequence $\langle \Phi_\alpha : \alpha < \kappa\rangle$ which uniformly witnesses these bijections exist. Also, is the condition that each $E_\alpha$ be an equivalence relation with all classes countable necessary in Question \ref{disjoint union smooth quotient isomorphic real product ordinal} and \ref{quotient by perfect tree isomorphism}? The following example of Holshouser-Jackson answers these questions.

\Begin{example}{some counterexample to previous}
Fix some recursive coding of binary relations on $\omega$ by reals. Let $\mathrm{WO}$ denote the collection of reals that code wellorderings on $\omega$. For $\alpha < \omega_1$, let $\mathrm{WO}_\alpha$ denote the reals coding wellorderings of ordertype $\alpha$. For $\alpha < \omega_1$, let $E_\alpha$ be the equivalence relation on $\reals$ defined by $x \ E_\alpha \ y$ if and only if $(x = y) \vee (x \notin \mathrm{WO}_\alpha \wedge y \notin \mathrm{WO}_\alpha)$. For each $\alpha < \omega_1$, $E_\alpha$ is $\borel$ bireducible to $=$. Hence $\reals \approx \reals \slash E_\alpha$. 

For each $\alpha < \omega_1$, if $x \in \mathrm{WO}_\alpha$, identify $[x]_{E_\alpha} = \{x\}$ with $x$. For each $\alpha < \omega_1$, $\reals \setminus \mathrm{WO}_\alpha$ is a single $E_\alpha$ equivalence class. Identify it with $\alpha$. Under this identification, one has a bijection of $\bigsqcup_{\alpha < \omega_1} \reals \slash E_\alpha$ with $\mathrm{WO} \sqcup \omega_1 \approx \reals \sqcup \omega_1$. 

$\reals \sqcup \omega_1$ is not in bijection with $\reals \times \omega_1$: Suppose $\Phi : \reals \sqcup \omega_1 \rightarrow \reals \times \omega_1$ is a bijection. $\pi_1[\Phi[\omega_1]]$ can be wellordered using $\Phi$ and the wellordering on $\omega_1$. ($\pi_1 : \reals \times \omega_1 \rightarrow \reals$ is the projection onto the first coordinate.) Under $\mathsf{AD}$, there is no uncountable sequence of distinct reals; hence, $\pi_1[\Phi[\omega_1]]$ is countable. Let $r \in \reals$ such that $r \notin \pi_1[\Phi[\omega_1]]$. $\Phi^{-1}[\{r\} \times \omega_1] \subseteq \reals$. But $\Phi^{-1}[\{r\} \times \omega_1]$ can be wellordered. This would give an uncountable sequence of distinct reals in $\reals$. Contradiction.

As mentioned above, if $\langle E_\alpha : \alpha < \omega_1\rangle$ was accompanied by a sequence of bijections $\langle \Phi_\alpha : \alpha < \omega_1\rangle$, then one can construction a bijection between $\bigsqcup_{\alpha < \omega_1} \reals \slash E_\alpha$ and $\reals \times \omega_1$. Thus, there cannot be such a sequence of bijections under $\mathsf{AD}$.

Note that $E_\alpha$ has exactly one uncountable class. This example shows Question \ref{disjoint union smooth quotient isomorphic real product ordinal} has a negative answer without the condition that each $E_\alpha$ has all countable classes. 

Let $p$ be a perfect tree such that $[p] \subseteq \reals \setminus \mathrm{WO}$. Then $\bigsqcup_{\alpha < \omega_1} [p] \slash E_\alpha \approx \omega_1$. $\omega_1$ is not in bijection with $\bigsqcup_{\alpha < \omega_1} \reals \slash E_\alpha \approx \reals \sqcup \omega_1$. Hence Question \ref{quotient by perfect tree isomorphism} has a negative answer if all the equivalence relations do not have all classes countable.
\end{example}

If all the equivalence relations in $\langle E_\alpha : \alpha < \kappa\rangle$ have all classes countable and $\reals \approx \reals \slash E_\alpha$, then $\bigsqcup_{\alpha < \kappa} \reals \slash E_\alpha$ contains a subset which is in bijection with $\reals \sqcup \omega_1$ but itself is not in bijection with $\reals \sqcup \omega_1$.

\Begin{fact}{R cup omega_1 and disjoint union smooth}
$(\mathsf{ZF + AD})$ Let $\kappa$ be an uncountable ordinal. Let $\langle E_\alpha : \alpha < \kappa\rangle$ be a sequence of equivalence relations on $\reals$ so that for each $\alpha < \kappa$, $E_\alpha$ has all classes countable and $\reals \approx \reals \slash E_\alpha$. Then $\reals \sqcup \kappa$ injects into $\bigsqcup_{\alpha < \kappa} \reals \slash E_\alpha$, but $\reals \sqcup \kappa$ is not in bijection with $\bigsqcup_{\alpha < \kappa} \reals \slash E_\alpha$. 
\end{fact}

\begin{proof}
Let $\bar{0} : \omega \rightarrow \{0\}$, be the constant $0$ function. For each $\alpha < \kappa$, identify $[\bar{0}]_{E_\alpha}$ with $\alpha$. Let $\Phi : \reals \rightarrow (\reals \slash E_0) \setminus [\bar{0}]_{E_0}$ be a bijection. Identify $\Phi(r)$ with $r$. Using this identification, there is a subset of $\bigsqcup_{\alpha < \kappa} \reals \slash E_\alpha$ which is in bijection with $\reals \sqcup \kappa$. 

Suppose there is a bijection $\Phi : \reals \sqcup \kappa \rightarrow \bigsqcup_{\alpha < \kappa} \reals \slash E_\alpha$. $\bigcup_{\alpha < \kappa} \Phi(\alpha)$ can be prewellordered by $x \sqsubseteq y$ if and only if the least $\alpha$ such that $x \in \Phi(\alpha)$ is less than or equal to the least $\alpha$ such that $y \in \Phi(\alpha)$. Each $\sqsubseteq$-class is countable. Fact \ref{no PWO uncountable set countable PWO class} implies that $\bigcup_{\alpha < \kappa} \Phi(\alpha)$ is countable. Let $r \in \reals$ with $r \notin \bigcup_{\alpha < \kappa} \Phi(\alpha)$. Let $X = \{[r]_{E_\alpha} : \alpha < \kappa\}$. $\Phi^{-1}[X]$ is an uncountable sequence of distinct reals in $\reals$. Contradiction.
\end{proof}

\cite{The-Cardinals-Below-CountableSubsetOmega} Theorem 2 shows that under $\mathsf{ZF + DC + AD_\reals}$, the only uncountable cardinals below $\reals \times \omega_1$ are $\omega_1$, $\reals$, $\reals \sqcup \omega_1$, and $\reals \times \omega_1$. Thus under these assumptions, if $\bigsqcup_{\alpha < \omega_1} \reals \slash E_\alpha$ is not in bijection with $\reals \times \omega_1$, then $\bigsqcup_{\alpha < \omega_1} \reals \slash E_\alpha$ cannot inject into $\reals \times \omega_1$. Moreover, \cite{The-Cardinals-Below-CountableSubsetOmega} Theorem 2 shows that under $\mathsf{ZF + DC + AD_\reals}$, the only uncountable cardinals below $[\omega_1]^\omega$ are $\omega_1$, $\reals$, $\reals \sqcup \omega_1$, $\reals \times \omega_1$, and $[\omega_1]^\omega$. An interesting question would be to compare the cardinality of $[\omega_1]^\omega$ and $\bigsqcup_{\alpha < \omega_1} \reals \slash E_\alpha$ when each $E_\alpha$ is an equivalence relation with all classes countable.

\Begin{fact}{no injection countable subsets omega1 into disjoint union quotient}
$(\mathsf{ZF + AD})$ Let $\langle E_\alpha : \alpha < \omega_1\rangle$ be a sequence of equivalence relations on $\reals$ such that each $E_\alpha$ has all classes $\bPi_1^0$. There is no injection of $[\omega_1]^\omega$ into $\bigsqcup_{\alpha < \omega_1} \reals \slash E_\alpha$. 
\end{fact}

\begin{proof}
Recall that $\mathcal{U}$ is Martin's cone measure on $\degrees$, the set of Turing degrees. For each $x \in \degrees$, let $\Lambda(x)$ denote the collection of countable $x$-admissible ordinals. For each $x \in \degrees$, let $\Gamma(x) \in [\omega_1]^\omega$ be the increasing sequence of the first $\omega$-many $x$-admissible ordinals.

Suppose $\Phi : [\omega_1]^\omega \rightarrow \bigsqcup_{\alpha < \omega_1} \reals \slash E_\alpha$ is an injection.

A sequence of Turing degrees $(x_n : n \in \omega)$ and a sequence $(\sigma_n : n \in \omega)$ in $\finBinarySequence$ will be constructed by recursion with the property that for all $n \in \omega$, $|\sigma_n | = n$, $\sigma_n \subset \sigma_{n + 1}$, and whenever $f \in [\Lambda(x_n)]^\omega$, there is some $r \in \Phi(f)$ so that $\sigma_n \subset r$. 

Let $x_0 = [\bar{0}]_T$, where $\bar{0}$ is the constant $0$ function. Let $\sigma_0 = \emptyset$. 

Suppose $x_n$ and $\sigma_n$ have been defined with the desired properties. Let $E^{n + 1}_0 = \{x \in \degrees : (\exists r \in \Phi(\Gamma(x)))(\sigma_n\hat{\ }0 \subseteq r)\}$ and $E^{n + 1}_1 = \{x \in \degrees : (\exists r \in \Phi(\Gamma(x)))(\sigma_n\hat{\ }1 \subseteq r)\}$. Note that the cone above $x_n$ is contained in $E^{n + 1}_0 \cup E^{n + 1}_1$. Since $\mathcal{U}$ is an ultrafilter, there is some $i \in 2$ and $x_{n + 1} \geq_T x_n$ so that $E^{n + 1}_i$ contains the cone above $x_{n + 1}$. Let $\sigma_{n + 1} = \sigma_n \hat{\ }i$, for this $i \in 2$. 

Let $f \in [\Lambda(x_{n + 1})]^\omega$. A result of Jensen (\cite{Admissible-Sets}) shows that for every increasing $\omega$-sequence of $x_{n + 1}$-admissible ordinals $f$, there is some $y \geq_T x_{n + 1}$ so that $\Gamma(y) = f$. Then $y \in E^{n + 1}_i$. Hence there is some $r \in \Phi(\Gamma(y)) = \Phi(f)$ so that $\sigma_{n + 1} \subseteq r$. 

Let $r = \bigcup_{n \in \omega} \sigma_n$. Let $z$ be the join $\bigoplus_{n \in \omega} x_n$. Suppose $f \in [\Lambda(z)]^\omega$. For all $n \in \omega$, there is some $r^f_n \in \Phi(f)$ so that $\sigma_n \subseteq r_n^f$. $\Phi(f)$ is an $E_\alpha$ class for some $\alpha < \omega$ so $\Phi(f)$ is $\bPi_1^0$. Since $r$ is the limit of $\{r^f_n : n \in \omega\} \subseteq \Phi(f)$, $r \in \Phi(f)$. It has been shown that for all $f \in [\Lambda(z)]^\omega$, $\Phi(f) \in \{[r]_{E_\alpha} : \alpha < \omega_1\} \approx \omega_1$. Then $\Phi$ induces an injection of $[\Lambda(z)]^\omega$ into $\omega_1$. This is impossible since such an injection would yield a wellordering of $\bbR$ since $\bbR$ injects into $[\Lambda(z)]^\omega$. 
\end{proof}

The above argument incorporates Martin's proof of the partition relation $\omega_1 \rightarrow (\omega_1)^\omega_2$. The following result captures the essential idea of the above argument.

\Begin{fact}{no injection countable sequence into disjoint union}
$(\mathsf{ZF + AD})$ Let $\kappa \in \mathrm{ON}$. Let $\langle E_\alpha : \alpha < \kappa\rangle$ be a sequence of equivalence relations on $\reals$ . Let $\Phi : [\omega_1]^\omega \rightarrow \bigsqcup_{\alpha < \kappa} \reals \slash E_\alpha$. Let $R \subseteq [\omega_1]^\omega \times \reals$ be defined by $R(f,x) \Leftrightarrow x \in \Phi(f)$. If $R$ has a uniformizing function then $\Phi$ is not an injection.
\end{fact}

\begin{proof}
Let $\Psi$ be a uniformizing function for $R$. 

For each $n \in \omega$, let $E^n_i = \{x \in \degrees : \Psi(\Gamma(x))(n) = i\}$. Since $\mathcal{U}$ is an ultrafilter, there is some $a_n \in 2$ such that $E_{a_n}^n \in \mathcal{U}$. Let $x_n \in \degrees$ be such that the cone above $x_n$ lies inside of $E^n_{a_n}$. Now suppose that $f \in [\Lambda(x_n)]^\omega$. A result of Jensen (\cite{Admissible-Sets}) states that for any such $f$, there is some $y \geq_T x_n$ so that $\Gamma(y) = f$. As $y \in E^n_{a_n}$, $\Psi(\Gamma(y))(n) = \Psi(f)(n) = a_n$. 

Let $r \in \reals$ be such that for all $n$, $r(n) = a_n$. Let $x = \bigoplus x_n$. If $f \in [\Lambda(x)]^n$, then $\Psi(f) = r$. 

It has been shown that there is an uncountable set $X \subseteq \omega_1$ and some real $r$ so that $\Psi[[X]^\omega] = \{r\}$. By definition of $R$, $\Phi[[X]^\omega] \subseteq \{[r]_{E_\alpha} : \alpha < \kappa\}$. The latter set is in bijection with $\kappa$. $[X]^\omega \approx [\omega_1]^\omega$. Therefore, $\Phi$ induces an injection of $[\omega_1]^\omega$ into the ordinal $\kappa$. As $\reals$ injects into $[\omega_1]^\omega$, this would imply that one could wellorder $\reals$. 
\end{proof}

Note that in Fact \ref{no injection countable sequence into disjoint union}, $R$ only needs to be uniformized on a set of cardinality $[\omega_1]^\omega$. To see this, suppose $R$ is uniformized on $Z \subseteq [\omega_1]^\omega$ of cardinality $[\omega_1]^\omega$. Let $L : [\omega_1]^\omega \rightarrow Z$ be a bijection. Let $\Phi' = \Phi \circ L$. The relation $R'$ associated to $\Phi'$ can be uniformized. Hence $\Phi'$ is not injective by Fact \ref{no injection countable sequence into disjoint union}. This implies $\Phi$ is not injective.

The class of equivalence relations with $\bPi_1^0$ classes is very restrictive. However, it does include equivalence relations with all finite classes. However, in such cases, there is a more natural argument: Fix some linear ordering $<$ of $\reals$. For $f \in [\omega_1]^\omega$, let $L(x)$ denote the $<$-least element of $\Phi(x)$ (which exists since $\Phi(x)$ is finite). Now apply Fact \ref{no injection countable sequence into disjoint union}.

\Begin{fact}{injection disjoint union into countable sequence}
(With Jackson.) Assume $\mathsf{ZF + AD^+}$. Let $\kappa \in \mathrm{ON}$ and $\langle E_\alpha : \alpha < \kappa\rangle$ be a sequence of equivalence relations on $\reals$ such that each $E_\alpha$ has all classes countable and $\reals \slash E_\alpha \approx \reals$. Then there is no injection $\Phi : [\omega_1]^\omega \rightarrow \bigsqcup_{\alpha < \kappa} \reals \slash E_\alpha$. 
\end{fact}

\begin{proof}
This is proved by verifying the uniformization condition of Fact \ref{no injection countable sequence into disjoint union}. Note that if $\langle E_\alpha : \alpha < \kappa\rangle$ is a sequence so that each $E_\alpha$ is an equivalence relation with all classes countable, then for any $\Phi$, the associated relation has all countable sections. 

Woodin's countable section uniformization states that every relation on $\reals \times \reals$ with countable section can be uniformized under $\mathsf{AD}^+$. In the present situation, the relations are on $[\omega_1]^\omega \times \reals$. Some modification of Woodin's ideas can be used to show countable section uniformization holds for such relations under $\mathsf{AD}^+$. The main ideas of Woodin's countable section uniformization on $\reals$ can be found in \cite{A-Trichotomy-Theorem-in-Natural} and \cite{Ramsey-Ultrafiler-and-Countable-to-One-Uniformation}. The details of this and other generalizations of countable section uniformization will appear elsewhere. 
\end{proof}

Originally, Theorem \ref{partial jonsson of dijoint union smooth quotient} was proved under $\mathsf{AD}^+$ using Woodin's countable section uniformization. However, it was observed that for the purpose of the J\'onsson property, one did not need total uniformization provided by Woodin's countable section uniformization but rather partial uniformization on a set of cardinality $\mathbb{R}$ (as by provided comeager uniformization) was adequate. As mentioned above, partial uniformization on a set of cardinality $[\omega_1]^\omega$ is adequate for the conclusion of Fact \ref{no injection countable sequence into disjoint union}. This suggests the following:

\Begin{question}{AD partial uniformization countable sequence}
Using just $\mathsf{AD}$, is it provable that for all relations $R \subseteq [\omega_1]^\omega \times \reals$ with countable sections, there is some $Z \subseteq [\omega_1]^\omega$ and $\Phi : Z \rightarrow \reals$ such that $|Z| = |[\omega_1]^\omega|$ and for all $z \in Z$, $R(z,\Phi(z))$?
\end{question}

The rest of this section will show the failure of the J\'onsson property for $(\reals \slash E_0) \times \kappa$ where $E_0$ is the equivalence relation from Definition \ref{some examples of equivalence relations} and $\kappa < \Theta$.

\Begin{fact}{subset R/E0 times ordinal has big projection}
$(\mathsf{ZF + AD})$ Suppose $A \subseteq (\reals \slash E_0) \times \kappa$ and $A \approx \reals \slash E_0$, where $\kappa$ is an ordinal. Let $\pi_1 : (\reals \slash E_0) \times \kappa \rightarrow \reals \slash E_0$ be the projection onto the first coordinate. Then $\pi_1[A] \approx \reals \slash E_0$. 
\end{fact}

\begin{proof}
Note that $A$ injects into $\pi_1[A] \times \kappa$. Hence $\reals \slash E_0$ injects into $\pi_1[A] \times \kappa$. Let $f : \reals \slash E_0 \rightarrow \pi_1[A] \times \kappa$ denote this injection. For each $\alpha < \kappa$, let $A_\alpha = \{x \in \reals : \pi_2(f([x]_{E_0})) = \alpha\}$, where $\pi_2 : (\reals \slash E_0) \times \kappa \rightarrow \kappa$ is the projection onto the second coordinate. Then $\bigcup_{\alpha < \kappa} A_\alpha = \reals$. By Fact \ref{union of meager sets}, there must be some $\alpha < \kappa$ so that $A_\alpha$ is nonmeager. Using the Baire property, $A_\alpha$ is comeager in some basic open set $O$. (Actually since $A_\alpha$ is $E_0$-invariant, it can be shown that $A_\alpha$ is comeager.) Hence $A_\alpha \supseteq \bigcap_{n \in \omega}D_n$, where $\langle D_n : n \in \omega\rangle$ is a sequence of topologically dense open sets relative to $O$. One can build an $E_0$-tree inside of $A_\alpha$. (See \cite{Definable-Combinatorics-Some-Borel-Equivalence-Relations} Definition 5.2.) This implies that there is a continuous reduction of $E_0$ into $E_0 \upharpoonright A_\alpha$. Hence $\reals \slash E_0$ injects into $A_\alpha \slash E_0$. Using $f$, $A_\alpha \slash E_0$ injects into $\pi_1[A] \times \{\alpha\} \approx \pi_1[A]$. It has been shown that $\reals \slash E_0$ injects into $\pi_1[A]$. Thus $\pi_1[A] \approx \reals \slash E_0$. 
\end{proof}

\Begin{fact}{R/E0 times ordinal is not 6 jonsson}
Let $\kappa < \Theta$. There is a $6$-J\'onsson function for $(\reals \slash E_0) \times \kappa$. 

$(\reals \slash E_0) \times \kappa$ is not J\'onsson.
\end{fact}

\begin{proof}
By Fact \ref{jonsson properties of some sets}, let $\Phi : [\reals \slash E_0]^3_= \rightarrow \reals \slash E_0$ be a $3$-J\'onsson map for $\reals \slash E_0$. Let $\Psi : \reals \rightarrow \kappa$ be a surjection. Since $=$ reduces into $E_0$, there is an injection $\Gamma : \reals \rightarrow \reals \slash E_0$. Let $\Lambda : [\reals \slash E_0]^3_= \rightarrow \kappa$ be defined by
$$\Lambda(x) = \begin{cases}
0 & \quad (\forall r \in \reals)(\Phi(x) \neq \Gamma(r)) \\
\Psi(r) & \quad \Phi(x) = \Gamma(r)
\end{cases}$$
Finally, let $\Upsilon : [(\reals \slash E_0) \times \kappa]^6_= \rightarrow (\reals \slash E_0) \times \kappa$ be defined by 
$$((x_1,\alpha_1), (x_2,\alpha_2), (x_3,\alpha_3), (x_4,\alpha_4), (x_5,\alpha_5), (x_6, \alpha_6)) \mapsto (\Phi(x_1,x_2,x_3), \Lambda(x_4,x_5,x_6))$$

Suppose $B \subseteq (\reals \slash E_0) \times \kappa$ is in bijection with $(\reals \slash E_0) \times \kappa$. Let $f : (\reals \slash E_0) \times \kappa \rightarrow B$ be a bijection. Let $A = f[(\reals \slash E_0) \times \{0\}]$. Then $A \approx \reals \slash E_0$. By Fact \ref{subset R/E0 times ordinal has big projection}, $\pi_1[A] \approx \reals \slash E_0$. 

Supppose that $(y,\beta) \in (\reals \slash E_0) \times \kappa$. Suppose $\Psi(r) = \beta$. Since $\Phi$ is a $3$-J\'onsson map and $\pi_1[A] \approx \reals \slash E_0$, one can find $((x_1,\alpha_1), (x_2,\alpha_2), (x_3,\alpha_3), (x_4,\alpha_4), (x_5,\alpha_5), (x_6, \alpha_6)) \in [A]^6_= \subseteq [B]^6_=$ so that $\Phi(x_1, x_2,x_3) = y$ and $\Phi(x_4,x_5,x_6) = \Gamma(r)$. Then $\Upsilon((x_1,\alpha_1), (x_2,\alpha_2), (x_3,\alpha_3), (x_4,\alpha_4), (x_5,\alpha_5), (x_6, \alpha_6)) = (y,\beta)$. $\Upsilon$ is a $6$-J\'onsson function for $(\reals \slash E_0) \times \kappa$. 
\end{proof}

\Begin{question}{jonsson function RE0 times ordinal questions}
\cite{Definable-Combinatorics-Some-Borel-Equivalence-Relations} showed that $\reals \slash E_0$ has no $2$-J\'onsson map but has a $3$-J\'onsson map. What is the least $n$ so that $(\reals \slash E_0)\times \kappa$ has a $n$-J\'onsson map, where $\kappa < \Theta$?

If $\kappa$ is any ordinal, is $(\reals \slash E_0) \times \kappa$ also not J\'onsson?
\end{question}

\bibliographystyle{amsplain}
\bibliography{references}

\end{document}